\documentclass[12pt,a4paper]{amsart} 
\usepackage{graphicx,latexsym,amssymb,color}
\theoremstyle{plain} 
\newtheorem{thm}{Theorem}[section] 
\newtheorem{lem}[thm]{Lemma} 
\newtheorem{prop}[thm]{Proposition} 
\newtheorem{cor}[thm]{Corollary} 
\theoremstyle{definition}
\newtheorem{defi}{Definition} 
\theoremstyle{remark} 
\newtheorem{rem}[thm]{Remark} 
 
\newtheorem{example}[thm]{Example}
\newtheorem*{acknowledgments}{Acknowledgments} 
\numberwithin{equation}{section}
\numberwithin{figure}{section}
\newcommand{\bd}{\begin{description}}   
\newcommand{\ed}{\end{description}} 
\newcommand{\ba}{\begin{array}}      \newcommand{\ea}{\end{array}} 
\newcommand{\bc}{\begin{center}}     \newcommand{\ec}{\end{center}} 
\newcommand{\be}{\begin{enumerate}}  \newcommand{\ee}{\end{enumerate}} 
\newcommand{\beq}{\begin{eqnarray}}  \newcommand{\eeq}{\end{eqnarray}} 
\newcommand{\beQ}{\begin{eqnarray*}} \newcommand{\eeQ}{\end{eqnarray*}} 
\newcommand{\bi}{\begin{itemize}}    \newcommand{\ei}{\end{itemize}}

\newcommand{\ov}{\overline}

\newcommand{\1}{ {{\mathbf 1}} }
\newcommand{\n}{ \{ 1,...,n \} }

\begin{document} 
\title[]{Self delta-equivalence for links whose 
Milnor's isotopy invariants vanish} 

\author[A. Yasuhara]{Akira Yasuhara} 
\address{Tokyo Gakugei University\\
         Department of Mathematics\\
         Koganeishi \\
         Tokyo 184-8501, Japan}
	 \email{yasuhara@u-gakugei.ac.jp}

\thanks{
The author is partially supported by a Grant-in-Aid for Scientific Research (C) 
($\#$18540071) of the Japan Society for the Promotion of Science.}

%
\begin{abstract} 
For an $n$-component link $L$, the Milnor's isotopy invariant  
is defined for each multi-index $I=i_1i_2...i_m~(i_j\in\n)$. 
Here $m$ is called the length. 
Let $r(I)$ denote the maximam number of times that any index appears. 
It is known that Milnor invariants with $r=1$ are link-homotopy 
invariant. N.~Habegger and X.~S.~Lin showed that two string links are a link-homotopc 
if and only if their Milnor invariants with $r=1$ coincide. 
This gives us that a link in $S^3$ is link-homotopic to a trivial link if and only if 
the all Milnor invariants of the link with $r=1$ vanish. 
Although Milnor invariants with $r=2$ are not link-homotopy invariants,   
T.~Fleming and the author showed that Milnor invariants with $r\leq 2$ 
are self $\Delta$-equivalence invariants. 
In this paper, we give a self $\Delta$-equivalence classification of 
the set of $n$-component links in $S^3$ whose Milnor invariants with length $\leq 2n-1$ 
and $r\leq 2$ vanish. As a corollary, we have that 
a link is self $\Delta$-equivalent to a trivial link if and only if 
the all Milnor invariants of the link with $r\leq 2$ vanish. 
This is a geometric characterization for links whose Milnor invariants 
with $\leq 2$ vanish. 
The chief ingredient in our proof is Habiro's clasper theory. 
We also give an alternate proof of a link-homotopy classification 
of string links by using clasper theory.   
\end{abstract} 

\maketitle 
\baselineskip=15pt

\bigskip
{\section{Introduction}}

For an $n$-component link $L$, {\em Milnor invariant}  
$\ov{\mu}_L(I)$ is defined for each multi-index $I=i_1i_2...i_m~(i_j\in\n)$ \cite{Milnor,Milnor2}.  
Here $m$ is called the {\em length} of $\ov{\mu}_L(I)$ and denoted by $|I|$. 
Let $r(I)$ denote the maximam number of times that any index appears. 
For example, $r(1123)=2,~r(1231223)=3$. It is known that if $r(I)=1$, 
then $\ov{\mu}_L(I)$ is a {\em link-homotopy} invariant \cite{Milnor}, where {link-homotopy} 
is an equivalence relation on links generated by self crossing changes. 
Similarly, for a string link $L$, Milnor invariant $\mu_L(I)$ is defined \cite{HL}. 
While Milnor invariants are not strong enough to give a link-homotopy classification for links, 
they are complete for  string links.  In fact, the following is known \cite{HL}. 

\begin{thm}[{\cite{HL}}]
\label{LHC} 
Two $n$-component string links $L$ and $L'$ are link-homotopic 
if and only if $\mu_L(I)=\mu_{L'}(I)$ for any $I$ with $r(I)=1$. 
\end{thm}

We will give an alternate proof in section~4 via 
clasper theory. Actually we will give representatives determined by 
Milnor link-homotopy invariants for the link-homotopy classes 
explicitely, see Theorem~\ref{SLC}. As a corollary, we have that for $n$-component string links 
$L$ and $L'$, and for a positive integer $k~(k\leq n)$, ${\mu}_L(I)={\mu}_{L'}(I)$ for 
any $I$ with $r(I)=1$ and $|I|\leq k$ 
if and only if $L$ and $L'$ are transformed into each other by combining link-homotopies and $C_k$-moves,  
see Corollary~\ref{SLC2}.

For a string link $L$, let $\mathrm{cl}(L)$ denote the {\em closure} of $L$. 
It follows from the definitions that $\mu_L(I)=\ov{\mu}_{\mathrm{cl}(L)}(I)$ if 
 $\mu_L(J)=0$ for any $J$ with $|J|<|I|$. 
 Since the Milnor invariants of trivial (string) links are $0$, 
this and Theorem~\ref{LHC} imply the following. 
The proposition below also follows from Milnor's link-homotopy 
classification theorem for Brunnian links \cite{Milnor}. 

\begin{prop}[{\cite[Section 5]{Milnor}}]
\label{VLH} 
A link $L$ in $S^3$ is link-homotopic to a trivial link if and only if 
$\ov{\mu}_L(I)=0$ for any $I$ with $r(I)=1$.  
\end{prop}

Although Milnor invariants with $r\geq 2$ are 
not necessarily link-homotopy invariants, they are generalized link-homotopy invariants. 
In fact, Fleming and the author \cite{FY} showed that Milnor invariants with $r\leq k$ are 
{\em self $C_k$-equivalence} invariants, where  
the (self) $C_k$-equivalence is an equivalence relation on (string) links 
generated by (self) $C_k$-moves defined as follows. 

\begin{figure}[!h]
\includegraphics[trim=0mm 0mm 0mm 0mm, width=.9\linewidth]
{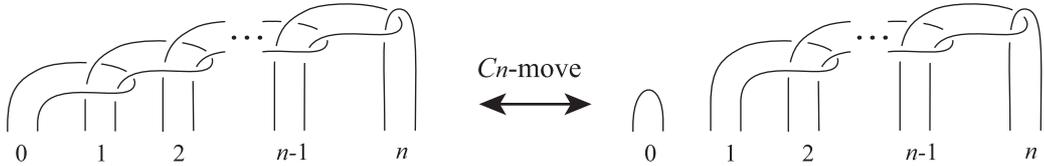}
\caption{A $C_n$-move involves $n+1$ strands of a link.   } \label{Cn}
\end{figure}

A {\em $C_n$-move}  is a local move on (string) links as illustrated in Figure~\ref{Cn}. 
 (A $C_1$-move is defined as the crossing change). 
These local moves were introduced by Habiro \cite{Hmt}.   
 A $C_n$-move is called a {\em self $C_n$-move} if the all strands in Figure~\ref{Cn} 
 belong to the same component of a (string) link \cite{SY1}. 

\begin{figure}[!h]
\includegraphics[trim=0mm 0mm 0mm 0mm, width=.45\linewidth]
{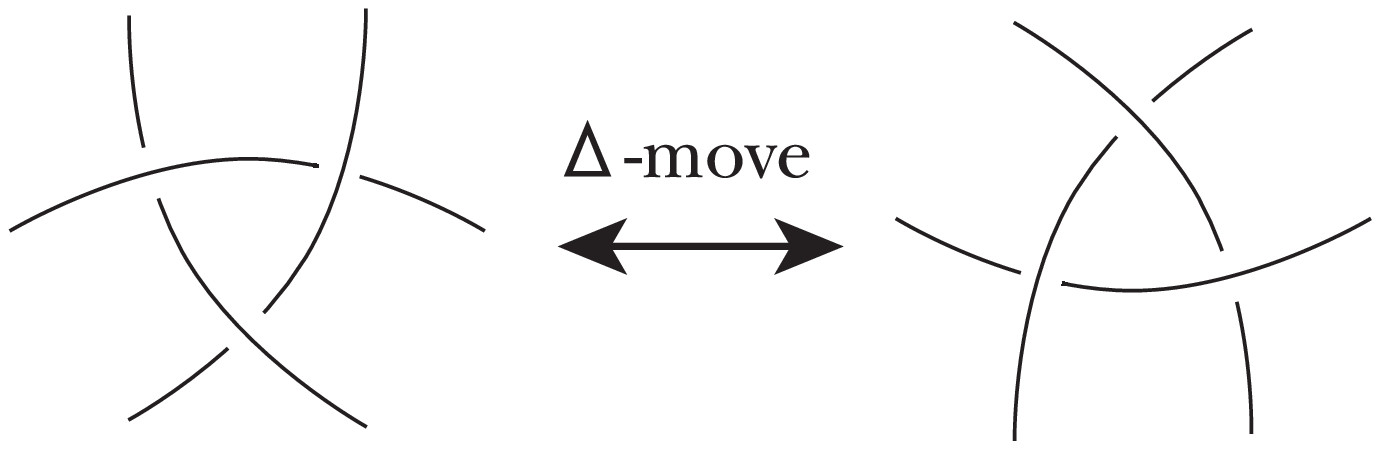}
\caption{}
\label{delta}
\end{figure}

The $C_n$-move (resp. self $C_n$-move) generates an equivalence relation on links, 
called the {\em $C_n$-equivalence} (resp. {\em self $C_n$-equivalence}).  
This notion can also be defined by using the theory of claspers (see section 2). 
The (self) $C_n$-equivalence relation becomes finer as $n$ increases, i.e., 
the  (self) $C_m$-equivalence implies the  (self) $C_k$-equivalence for $m>k$.  
We remark that (self) $C_2$-move is same as (self) {\em $\Delta$-move} defined by \cite{MN}. 
The $\Delta$-move is defined as a local move as illustrated in 
Figure~\ref{delta}. 
We call the (self) $C_2$-equivalence the {\em (self) $\Delta$-equivalence}.

A self $\Delta$-equivalence classification of {\em 2-component} links 
was shown by Y. Nakanishi and Y. Ohyama \cite{NO}. 
It is still open for links with at least 3 components. 
Here we give the following theorem. 

\begin{thm}\label{SDC} 
Let $L$ and $L'$ be $n$-component links. Suppose that $\ov{\mu}_L(I)=\ov{\mu}_{L'}(I)=0$ 
for any $I$ with $|I|\leq 2n-1$ and $r(I)\leq 2$. 
Then $L$ and $L'$ are self $\Delta$-equivalent if and only if 
$\ov{\mu}_L(J)=\ov{\mu}_{L'}(J)$ for any $J$ with $|J|= 2n$ and $r(J)= 2$. 
\end{thm}

\begin{rem}
(1)~The \lq only if' part follows directly from 
the fact that Milnor invariants with $r\leq k$ are self $C_k$-equivalence 
invariants \cite{FY}.  \\
(2)~In the last section, we characterize $n$-component links whose Milnor 
invariants of length $\leq 2n-1$ and $r\leq 2$ vanish. 
More precisely, the Milnor invariants of an $n$-component link  
with length $\leq 2n-1$ and $r\leq 2$ vanish if and only if, for any integer $i$ in $\{1,...,n\}$, 
it is self $\Delta$-equivalent to 
a Brunnian link $L_i$ such that the $i$th component $K$ of $L_i$ is null-homotopic in 
$S^3\setminus(L_i-K)$ (Theorem~\ref{char}). 
As an example, we will give a 3-component Brunnian link $L=K_1\cup K_2\cup K_3$ 
such that $K_1$ is not null-homotopic in $S^3\setminus(L-K_1)$ 
and $K_i$ is null-homotopic in $S^3\setminus(L-K_i)~(i=2,3)$ (Example~\ref{example}). 
In particular, $L$ is link-homotopic to a trivial link. 
There is no such a link with 2 components, i.e., if a 2-component link is link-homotopic to 
a trivial link, then it is self $\Delta$-equivalent to a Brunnian link $K_1\cup K_2$ such that 
$K_i$ is null-homotopic in $S^3\setminus K_j~(\{i,j\}=\{1,2\})$. 
\end{rem}

For 2-component links, 
Proposition~\ref{VLH} can be generalized \cite{NO}. 
Theorem~\ref{SDC} gives us  the following corollary 
which is a generalization of Proposition~\ref{VLH} for links with 
arbitrarily many components. 
This gives us a geometric characterization for links whose Milnor invariants 
with $\leq 2$ vanish.

\begin{cor}\label{COR1} 
A link $L$ is self $\Delta$-equivalent to a trivial link if and only if 
$\ov{\mu}_L(I)=0$ for any $I$ with $r(I)\leq 2$. 
\end{cor}

 \begin{rem}
(1)~This corollary gives an affirmative answer for an open question 
remained in \cite{FY}.\\
(2)~For string links, Corollary\ref{COR1} does not hold, i.e., there are 
 2-string links such that their Milnor invariants $\mu(I)$ with 
 $r(I)\leq 2$ vanish and they are not self $\Delta$-equivalent to 
 a trivial string link \cite{FY2}.\\
(3)~Since $C_k$-move ($k\geq 3$) is not unknotting operation, 
 it is impossible to generalize the corollary above. 
 It is reasonable to consider the following question: 
 If $\ov{\mu}_L(I)=0$ for any $I$ with $r(I)\leq k$, then 
 is $L$ self $C_k$-equivalent to a completely split link? 
 Fleming and the author gave a negative answer to the question \cite{FY}. 
 In fact, there is a 2-component boundary link $L$ such that 
 $L$ is not self $C_3$-equivalent to a split link. 
 Note that the all Milnor invariants of a boundary link vanish. 
 \end{rem}

By combining Lemma~\ref{parallel} (\cite[Theorem 7]{Milnor2}), Proposition~\ref{VLH} 
and Corollary~\ref{COR1}, we have the following corollary. 

\begin{cor}\label{COR2} 
Let $L$ be an $n$-component link and let $L(2)$ be 
a $2n$-component link obtained from $L$ by 
replaceing each component of $L$ with zero framed 2 parallel 
copies of it. 
Then $L$ is self $\Delta$-equivalent to a trivial link if and only if 
$L(2)$ is link-homotopic to a trivial link. 
\end{cor}
 
\begin{rem}
For an $n$-component link, let $L(k)$ be a 
$kn$-component link obtained from $L$ by 
replacing each component of $L$ with zero framed $k$ parallel 
copies of it.  In the proof of Theorem 2.1 in \cite{FY}, 
it is shown that if two links $L$ and $L'$ are self $C_k$-equivalent, then 
$L(k)$ and $L'(k)$ are link homotopic. 
So one might expect that if $L(2)$ and $L'(2)$ are link homotopic, 
then $L$ and $L'$ are self $\Delta$-equivalent. 
But this is not true. The reason is the follwing: 
There are 2-component links $L$ and $L'$ such that 
they are concordant and are not self $\Delta$-equivalent \cite{NS}, \cite{NSY}. 
The fact that $L$ and $L'$ are concordant implies that 
$L(2)$ and $L'(2)$ are concordant.   
Since link-concordance implies link-homotopy \cite{Gif}, \cite{Gol},  
$L(2)$ and $L'(2)$ are link-homotopic. 
\end{rem}

An $n$-component link $L=K_1\cup\cdots\cup K_n$ is called a 
{\em boundary link} if there is a disjoint union 
$X=F_1\cup \cdots\cup F_n$ of orientable surfaces such that 
$\partial X=L$ and $\partial F_i=K_i\ (i=1,2,...,n)$. 
An $n$-component link $L$ is called a 
{\em homology boundary link}  if $\pi_1(S^3\setminus L)$ 
admits an epimorphism from $\pi_1(S^3\setminus L)$ to a free group 
of rank $n$ \cite{Smy}. An every boundary link is a homology boundary link. 
T.~Shibuya and the author showed that all boundary links are 
self $\Delta$-equivalent to trivial links \cite{SY2}. 
In \cite{Shi}, Shibuya showed that all ribbon links are self $\Delta$-equivalent 
to trivial links. 

Whether the homology boundary links are self $\Delta$-equivalent to trivial links 
and whether the slice links are self $\Delta$-equivalent to trivial links
have remained as open questions. 
Since all Milnor invariants of homology boundary links vanish, and since 
Milnor invariants are concordance invariants, 
we have the following corollary, which give affirmative answers for the open questions. 
 
\begin{cor}\label{COR3} 
If $L$ is concordant to a homology boundary link, then $L$ is  
self $\Delta$-equivalent to a trivial link. 
\end{cor}


\begin{acknowledgments}
The author would like to thank Jean-Baptiste Meilhan for many useful discussions. 
The first joint work \cite{MeiY} leads to this work. 
He is also very grateful to Professor Tim Cochran for helpful comments.  
\end{acknowledgments}

\bigskip
{\section{Clasper} }
 
Let us briefly recall from \cite{H} the basic notions of clasper theory for (string) links.  
In this paper, we essentially only need the notion of $C_k$-tree.  
For a general definition of claspers, we refer the reader to \cite{H}.  
\begin{defi}\label{defclasp}
Let $L$ be a link in $S^3$ (resp. a string link in $D^2\times I$).  
An embedded disk $F$ in $S^3$ (resp. $D^2\times I$) 
is called a {\em tree clasper} for $L$ if 
it satisfies the following (1), (2) and (3):\\
(1) $F$ is decomposed into disks and bands, called {\em edges}, each of which 
connects two distinct disks.\\
(2) The disks have either 1 or 3 incident edges, called {\em leaves} or 
{\em nodes} respectively.\\
(3) $L$ intersects $F$ transversely and the intersections are contained 
in the union of the interior of the leaves. \\
The \emph{degree} of a tree clasper is the number of the leaves \emph{minus} $1$.  
(In \cite{H}, a tree clasper and a leaf are called 
a {\em strict tree clasper} and a {\em disk-leaf} respectively.) 
A degree $k$ tree clasper is called a {\em $C_k$-tree} (or a {\em $C_k$-clasper}).  
A $C_k$-tree is \emph{simple} if each leaf intersects $L$ at one point.  
\end{defi}

We will make use of the drawing convention for claspers of \cite[Fig. 7]{H}, except 
for the following: $\oplus$ (resp. $\ominus$) on an edge represents a positive 
(resp. negative) half-twist. (This replaces the 
convention of a circled $S$ (resp. $S^{-1}$) used in \cite{H}).    

\begin{figure}[!h]
\includegraphics[trim=0mm 0mm 0mm 0mm, width=.8\linewidth]
{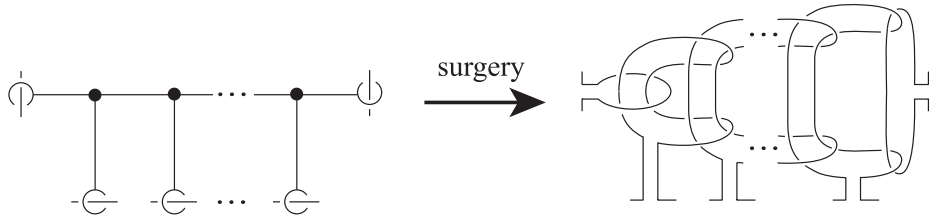}
\caption{Surgery along a simple $C_k$-tree.} \label{milnor-tangle}
\end{figure}

Given a $C_k$-tree $T$ for a link $L$ in $S^3$, there is a procedure to construct  a 
framed link $\gamma(T)$ in a regular neighborhood of $T$. 
\emph{Surgery along $T$} means surgery along $\gamma(T)$.  
Since there exists a canonical homeomorphism between $S^3$ and the manifold $S^3_{\gamma(T)}$,    
surgery along the $C_k$-tree $T$ can be regarded as a local move on $L$ in $S^3$.  
We say that the resulting link $L_T$ in $S^3$ is {\em obtained by surgery along $T$}.  
In particular, surgery along a simple $C_k$-tree illustrated in Figure~\ref{milnor-tangle} 
is equivalent to band-summing a copy of the $(k+1)$-component 
Milnor link (see \cite[Fig. 7]{Milnor}), 
and is equivalent to a $C_k$-move as defined in the introduction (Figure~\ref{Cn}).  
Similarly, for a disjoint union of trees $T_1\cup \cdots \cup T_m$,  we can define 
$L_{T_1\cup\cdots\cup T_m}$ as a link by surgery along $T_1\cup\cdots\cup T_m$.
A $C_k$-tree $T$ having the shape of the tree clasper in 
Figure~\ref{milnor-tangle} is called \emph{linear}, 
and the left-most and right-most leaves of $T$ in Figure 2.1 are called \emph{ends} of $T$.  
Ends of $T$ are not uniquely determined. There are 4 choices for an each linear tree. 

It is known that the $C_k$-equivalence as defined in section 1 coincides with 
the equivalence relation on links generated 
by surgery along $C_k$-trees and ambient isotopies.  
Two (string) links $L$ and $L'$ are $C_k$-equivalent if and only if 
there is a disjoint union of simple $C_k$-trees $G_1\cup\cdots\cup G_m$ such that 
$L'$ is ambient isotopic to $L_{G_1\cup\cdots\cup G_m}$ \cite[Theorem 3.17]{H}. 

\begin{defi} \label{r9}
  Let $L=K_1\cup\cdots\cup K_n$ be an $n$-component (string) link.  
A (simple) $C_k$-tree $T$ for $L$ is a \emph{(simple) $C^a_k$-tree} 
(resp. {\em $C_k^d$-tree, $C_k^s$-tree}) if it satisfies the following: \\
(1)~For each disk-leaf $f$ of $T$, $f\cap L$ is contained in a single component of $L$, and \\
(2)~ $|\{i~|~T\cap K_i\neq\emptyset\}|=n$ (resp. $=k+1,~1$).\\
Note that $n$ is the number of the components of $L$ and that 
$k+1$ is the number of leaves of $T$. 
If $T$ is simple, $T$ always satisfies the condition (1). 
The \emph{$C^*_k$-equivalence} ($*=a,d,s$) is an equivalence relation on (string) links 
generated by surgery along $C^*_k$-trees 
and ambient isotopies.  Note that 
$C_k^s$-equivalence is same as self $C_k$-equivalence. 
For a simple $C_k$-tree $T$, the set $\{i~|~T\cap K_i\neq\emptyset\}$ is called 
{\em index} of $T$, and denote it by $\mathrm{index}(T)$. 
And let $r_i(T)$ be the number of intersection points in $T\cap K_i$ $(i=1,...,n)$.  
The $(C_l^s+C_k^*)$-equivalence $(C_k^*=C_k,~C_k^a,$ or $C_k^d)$  
is an equivalence relation on (string) links generated by surgery along 
$C_l^s$- or $C_k^*$-trees. 
By the arguments similar to that in the proof of \cite[Theorem 3.17]{H}, we have 
that two (string) links $L$ and $L'$ are $(C_l^s+C_k^*)$-equivalent if and only if 
there is a disjoint union of simple $C_l^s$- or $C_k^*$-trees 
$T_1\cup\cdots\cup T_m$ such that 
$L'$ is ambient isotopic to $L_{T_1\cup\cdots\cup T_m}$. 
We use the notation $L \stackrel{C^*_{k}}{\sim} L'$ 
(resp. $L \stackrel{C_l^s+C^*_{k}}{\sim} L'$) for $C^*_k$-equivalent 
(resp. $(C_l^s+C_k^*)$-equivalent) links $L$ and $L'$.  
\end{defi}

Recall that a string link is a tangle without closed components 
(see \cite{HL} for a precise definition).  
The set of ambient isotopy classes of the $n$-component string links 
has a monoid structure with composition given by the \emph{stacking product}, 
denoted by $*$, and 
with the trivial $n$-component string link $\1_n$ as unit element. 

In the following, we give several lemmas. 
The proofs of Lemmas~\ref{cc}, \ref{slide} and 
\ref{twist} 
are essentially given in \cite{H} (see also section 1.4 in \cite{these}), and 
Lemma~\ref{ihx} essentially shown in \cite{G} (see also \cite{CT}, \cite{meilhan}), 
while they did not care about $r_j$ of claspers in \cite{CT}, \cite{G},  
\cite{H}, \cite{these}, \cite{meilhan}.   
If we follow their proofs with paying attention to $r_j$, 
we will see the proof of Lemmas~\ref{cc}, \ref{slide}, 
\ref{twist} and \ref{ihx}.

\begin{lem}[cf. {\cite[Propositions 4.5, 4.6]{H}}]\label{cc} 
Let $T$ be a  simple $C_k$-tree for an $n$-component  
(string) link $L$, and let $T'$ (resp. $T''$, and $T'''$) 
be obtained from $T$ by changing   
a crossing of an edge and the $i$th component $K_i$ of $L$ (resp. an edge of $T$, and
an edge of another simple clasper $G$) (see Figure~\ref{crossingchange}).  
Then \\
(1)~$L_T \stackrel{C_{k+1}}{\sim} L_{T'}$, and the $C_{k+1}$-equivalence 
is realized by surgery along 
simple $C_{k+1}$-trees with indices $\mathrm{index}(T)\cup\{i\}$ and  
$r_j\geq r_j(T)~(j=1,...,n)$. \\ 
(2)~$L_T \stackrel{C_{k+1}}{\sim} L_{T''}$, 
$L_{T\cup G} \stackrel{C_{k+1}}{\sim} L_{T'''\cup G}$, 
and the $C_{k+1}$-equivalence 
is realized by surgery along 
simple $C_{k+1}$-trees with $r_j\geq r_j(T)~(j=1,...,n)$.
\begin{figure}[!h]
\includegraphics[trim=0mm 0mm 0mm 0mm, width=.9\linewidth]
{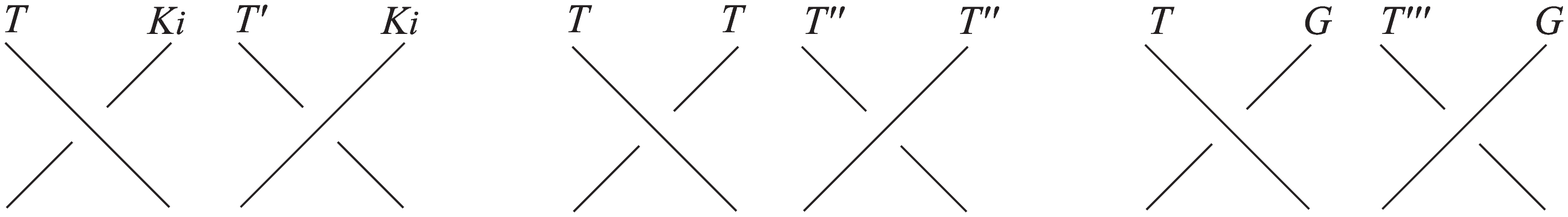}
\caption{} \label{crossingchange}
\end{figure}
\end{lem}

\begin{lem}[cf. {\cite[Propositions 4.4]{H}}]\label{slide} 
Let $T_1$ (resp. $T_2$) be a simple $C_k$-tree (resp. $C_l$-tree) for an 
$n$-component (string) link $L$, 
and let $T_1'$ be obtained from $T_1$ by 
sliding a leaf of $T_1$ over a leaf of $T_2$ (see Figure~\ref{sliding}).  
Suppose that $k\geq l$. 
Then $L_{T_1\cup T_2} \stackrel{C_{k+1}}{\sim} L_{T'_1\cup T_2}$, and 
the $C_{k+1}$-equivalence is realized by surgery along 
simple $C_{k+1}$-trees with $r_j\geq r_j(T_1)~(j=1,...,n)$.    

\begin{figure}[!h]
\includegraphics[trim=0mm 0mm 0mm 0mm, width=.4\linewidth]
{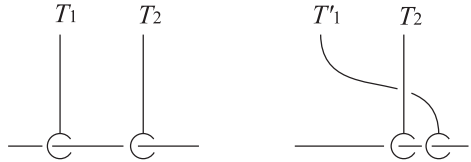}
  \caption{Sliding a leaf over another leaf.  }\label{sliding}
 \end{figure}
\end{lem}

\begin{lem}[cf. {\cite[Claim in p-36]{H}}]\label{twist} 
 Let $T$ be a simple $C_k$-tree for $\1_n$ and let $\ov{T}$ be a simple $C_k$-trees 
obtained from $T$ by adding a half-twist on an edge. Then 
 $(\1_n)_T* (\1_n)_{\ov{T}} \stackrel{C_{k+1}}{\sim} \1_n$, and  
the $C_{k+1}$-equivalence is realized by surgery along simple $C_{k+1}$-trees 
with $r_j\geq r_j(T)~(j=1,...,n)$.
\end{lem}


\begin{lem}[cf. {\cite[Theorem 6.7]{G}}, {\cite[Lemma 2.9]{meilhan}}] \label{ihx} 
Consider simple $C_k$-trees $T_I$, $T_H$ and $T_X$ for $\1_n$ 
which differ only in a small ball 
as illustrated in Figure~\ref{ihxmove}. Then 
$(\1_n)_{T_I} \stackrel{C_{k+1}}{\sim} (\1_n)_{T_H}* (\1_n)_{T_X}$, and 
the $C_{k+1}$-equivalence is realized by surgery along simple $C_{k+1}$-trees 
with $r_j\geq r_j(T_I)~(j=1,...,n)$.

\begin{figure}[!h]
\includegraphics[trim=0mm 0mm 0mm 0mm, width=.45\linewidth]
{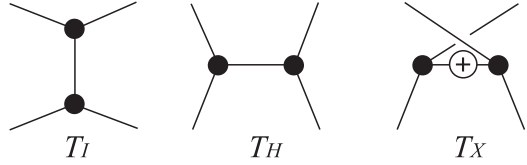}
  \caption{The IHX relation for $C_k$-trees.  }\label{ihxmove}
 \end{figure}
\end{lem} 

By combining the proof of \cite[Claim in p-26]{H} and 
\cite[Propositions 4.4, 4.5 and 4.6]{H}, we have the following. 

\begin{figure}[!h]
\includegraphics[trim=0mm 0mm 0mm 0mm, width=.4\linewidth]
{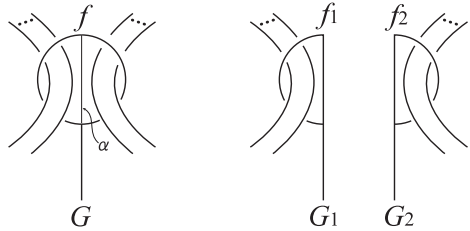}
  \caption{Splitting a leaf.  }\label{splitting}
 \end{figure}

\begin{lem}[cf. {\cite[Claim in p-26]{H}}]\label{split} 
Let $G$ be a $C_k$-tree for $\1_n$. Let $f_1$ and $f_2$ be two disks obtained by splitting a leaf $f$ of $G$ 
along an arc $\alpha$ as shown in Figure~\ref{splitting} 
(i.e., $f=f_1\cup f_2$ and $f_1\cap f_2=\alpha$). 
Then, $(\1_n)_G \stackrel{C_{k+1}}{\sim} (\1_n)_{G_1}* (\1_n)_{G_2}$, 
where $G_i$ denotes the $C_k$-tree for $\1_n$ obtained from $G$ by replacing $f$ with $f_i$ ($i=1,2$).  

\end{lem}

An $n$-component (string) link $L$ is {\em Brunnian} if every proper 
sublink of $L$ is trivial.  
In particular, any trivial (string) link is Brunnian.   
The $n$-component Brunnian (string) links are characterized  by $C_{n-1}^a$-equivalence 
as follows.

\begin{prop}[\cite{Hb,MY}] \label{brunnian} 
 Let $L$ be an $n$-component (string) link in $S^3$.  Then $L$ is Brunnian if and only if 
 $L$ is obtained from a trivial (string) link by surgery along simple $C_{n-1}^a$-trees.  
\end{prop}

By the argument similar to that in the proof of Theorem 1.2 in \cite{MY}, 
we have the following lemma. 
In \cite{MY}, they proved it with using \lq band description' defined 
by K.~Taniyama and the author \cite{TY2}. 
Here we give a proof with using clasper. 

\begin{lem}[cf. {\cite[Theorem 1.2]{MY}}] \label{flhb} 
Let $L$ be an $n$-component Brunnian link in $S^3$. 
If $L$ is obtained from a trivial link $O$ by surgery along 
$C_1^s$-trees with indices $\{i\}$, 
then $L$ is obtained from $O$ by surgery along simple $C_n^a$-trees with $r_i= 2$. 
\end{lem}

\begin{proof}
Set $O=O_1\cup \cdots\cup O_n$. 
It is enough to consider the case when $i=1$.
There is a disjoint union $F_1$ of simple $C_1^s$-trees with indices $\{1\}$ such that 
$L=O_{F_1}$. 
Note that $r_1=2$ for all $C_1^s$-trees in $F_1$.
 
Since $L$ is Brunnian, $L\setminus O_2$ is trivial. 
This implies that a split sum of $L\setminus O_2$ and $O_2$ is trivial.  
Hence $L$ can be deformed into trivial by crossing changes between 
$O_2$ and edges of $C_1^s$-trees of $F_1$. 
By Lemma~\ref{cc}, we have that $L$ is obtained from $O$ by 
surgery along a disjoint union $F_2$ of simple $C_2$-trees with indices $\{1,2\}$ and $r_1=2$. 
So we have $L=O_{F_2}$. 

Since $L\setminus O_3$ is trivial,  
$L$ can be deformed into trivial by crossing changes between 
$O_3$ and edges of $C_2$-trees in $F_2$. 
By Lemma~\ref{cc}, there is a disjoint union $F_3$ of simple $C_3$-trees with indecis 
$\{1,2,3\}$ and $r_1=2$ such that $L=O_{F_3}$. 

Repeating this step, we have that there is a disjoint union $F_n$ of 
simple $C_n$-trees with indices $\{1,...,n\}$ and $r_1=2$ such that 
$L=O_{F_n}$. This completes the proof. 
\end{proof}

By the arguments similar to that in the proof of \cite[Proposition 3.1]{FY}, 
we have

\begin{prop}[cf. {\cite[Proposition 3.1]{FY}}]
\label{ckindex} 
A (string) link $L'$ is obtained from $L$ by surgery along 
a simple $C_l$-tree with leaves $f_1,f_2,...,f_{l+1}$, then 
for any $k~(1\leq k<l)$ and 
any subset $\{w_1,...,w_{k+1}\}\subset \{1,...,{l+1}\}$, 
there are simple $C_k$-trees $T_j~(j=1,...,m)$ with leaves 
$f_{j1},...,f_{j(k+1)}$ such that $f_{ji}$ and $f_{w_i}$ 
grasp the same component of $L$  
for each $i(=1,...,k+1)$, and that $L'$ is obtained from 
$L$ by surgery along $T_1,...,T_m$.   
\end{prop} 

A simple $C_k$-tree $T$ is a {\em $C_k^{(l)}$-tree} if 
$\max\{r_j(T)~|~j=1,...,n\}\geq l$. 
Two links $L$ and $L'$ are {\em $C_k^{(l)}$-equivalent} if $L$ is 
obtained from $L'$ by ambient isotopy and surgery along 
simple $C_k^{(l)}$-trees. 
The following proposition is a corollary of Proposition~\ref{ckindex}.

\begin{prop}[cf. {\cite[Proposition 3.1]{FY}}]\label{sck} 
If two (string) links are $C_n^{(k+1)}$-equivalent, then 
they are self $C_k$-equivalent.   
Moreover, for some $i$, if the $C_n^{(k+1)}$-equivalence 
is realized by surgery along simple $C_n^{(k+1)}$-trees with 
$r_i\geq k+1$, then the self $C_k$-equivalence is realized 
by surgery along simple $C_k^s$-trees with indices $\{i\}$. 
\end{prop} 

\bigskip
{\section{Milnor invariants} }

J. Milnor defined in \cite{Milnor} a family of invariants of oriented, ordered 
links in $S^3$, known as Milnor's $\ov{\mu}$-invariants. 

Given an $n$-component link $L$ in $S^3$, denote by $G$ the fundamental group of 
$S^3\setminus L$, and by $G_q$ the $q\mathrm{th}$ subgroup of the lower central series of $G$.  
We have a presentation  of $G/ G_q$ with $n$ generators, given by a meridian 
$m_i$ of the $i\mathrm{th}$ component of $L$.  
So for $1\le i\le n$, the longitude $l_i$ of the $i\mathrm{th}$ component of $L$ is expressed 
modulo $G_q$ as a word  
in the $m_i$'s (abusing notations, we still denote this word by $l_i$).  

The \emph{Magnus expansion} $E(l_i)$ of $l_i$ is the formal power series in 
non-commuting variables $X_1,...,X_n$ obtained by 
substituting $1+X_j$ for $m_j$ and $1-X_j+X_j^2-X_j^3+\cdots$ for $m_j^{-1}$, $1\le j\le n$.  

Let $I=i_1 i_2 ...i_{k-1} j~(k\leq q)$ be a multi-index (i.e., a sequence of possibly repeating 
indices) among $\n$. 
Denote by $\mu_L(I)$ the coefficient of $X_{i_1}\cdots X_{i_{k-1}}$ in the Magnus expansion $E(l_j)$.  
\emph{Milnor invariant} $\ov{\mu}_L(I)$ is the residue class of $\mu_L(I)$ modulo the greatest common divisor of 
all Milnor invariants $\mu_L(J)$ such that $J$ is obtained from $I$ by removing at least one index.  
As we mentioned in section 1, 
$|I|=k$ is called the {length} of Milnor invariant $\ov{\mu}_L(I)$.  
 
The indeterminacy comes from the choice of the meridians $m_i$.  
Equivalently, it comes from the indeterminacy of 
representing the link as the closure of a string link \cite{HL}.  
Indeed, $\mu(I)$ is a well-defined invariant for string links.  

The following 4 lemmas play an important roles in calculating 
Milnor invariants. 

\begin{figure}[!h]
\includegraphics[trim=0mm 0mm 0mm 0mm, width=.4\linewidth]
{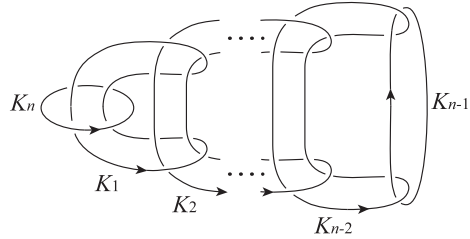}
\caption{Milnor link} \label{milnorlink1}
\end{figure}

\begin{lem}[{\cite[section 5]{Milnor}}] \label{Milnorlink} 
Let $M_n=K_1\cup\cdots\cup K_n$ be the $n$-component Milnor link 
as illustrated in Figure~\ref{milnorlink1}. Then the Milnor invariants of 
length $\leq n-1$ vanish, and 
\[\ov{\mu}_{M_n}(i_1i_2...i_{n-2}~n-1~n)=
\left\{\begin{array}{ll}
1& \text{ if $i_1i_2... i_{n-2}=12... n-2$},\\
0& \text{ otherwise}.
\end{array}
\right.\]
\end{lem}

\begin{lem}[{\cite[Theorem 7]{Milnor2}}]\label{parallel} 
Let $L'$ be a link obtained from a link $L$ by taking 
the appropriate number of zero framed parallels of the components of $L$. 
Suppose the $i$th component of $L'$ corresponds to the $h(i)$th component 
of $L$, then
\[ \ov{\mu}_{L'}(i_1i_2... i_m)=
\ov{\mu}_{L}(h(i_1)h(i_2)... h(i_m)).\]
\end{lem}

\begin{lem}[{\cite[Lemma 3.3]{MeiY}}]\label{add} 
Let $L$ and $L'$ be $n$-component string links such that all Milnor invariants of $L$ 
(resp. $L'$) of length $\le m$ 
(resp. $\le m'$) vanish.  
Then $\mu_{L* L'}(I)=\mu_{L}(I)+\mu_{L'}(I)$ for all $I$ of length $\leq m+m'$.  
\end{lem}

\begin{lem}[{\cite[Theorem 7.2]{H}}] \label{milnCk} 
The Milnor invariants of length $\leq k$ for (string) links are invariants of 
the $C_{k}$-equivalence.  
\end{lem}

\bigskip
{\section{Link-homotopy of string links}}

\begin{figure}[!h]
\includegraphics[trim=0mm 0mm 0mm 0mm, width=.8\linewidth]
{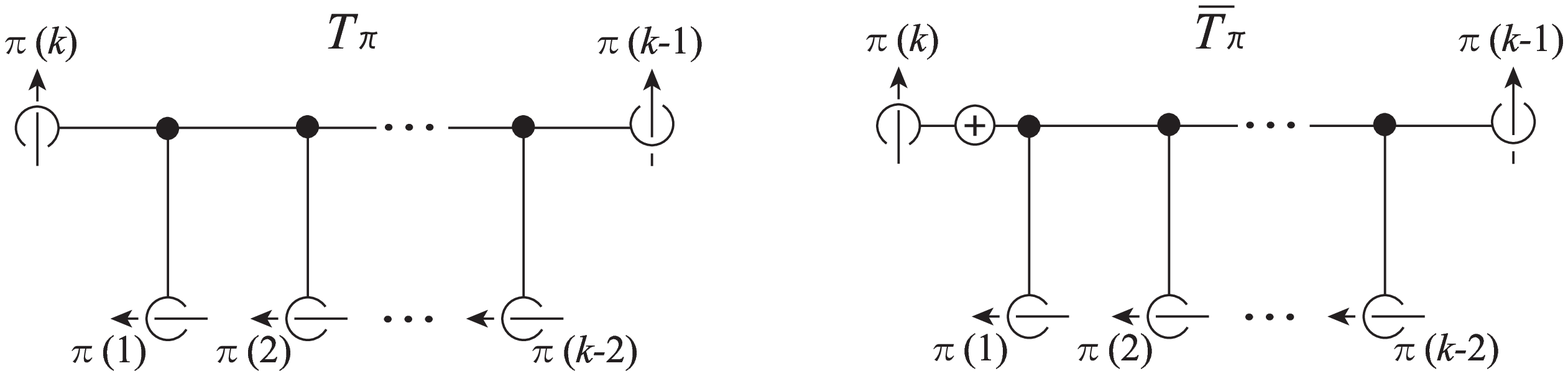}
\caption{${\pi(i)}$ means $\pi(i)$th component of $\1_n$} \label{baseT}
\end{figure}

\begin{figure}[!h]
\includegraphics[trim=0mm 0mm 0mm 0mm, width=.45\linewidth]
{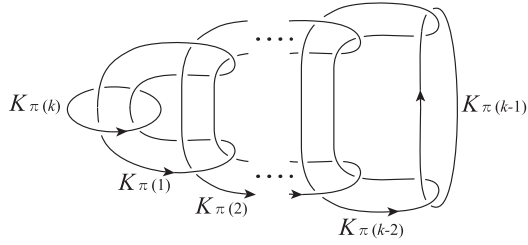}
\caption{Milnor link $M_{\pi}=K_{\pi(1)}\cup\cdots\cup K_{\pi(k)}$} \label{milnorlink2}
\end{figure}

Let $\pi:\{1,...,k\}\longrightarrow\n~(k\leq n)$ be an injection such that
$\pi(i)<\pi(k-1)<\pi(k)~(i\in\{1,...,k-2\})$, and let $\mathcal{F}_k$ be the 
set of such injections. 
 For $\pi\in\mathcal{F}_k$, let $T_{\pi}$ and $\ov{T}_{\pi}$ be simple 
 $C_{k-1}^d$-trees as illustrated in Figure~\ref{baseT}, and 
 set $V_{\pi}=(\1_n)_{T_{\pi}}$ and $V^{-1}_{\pi}=(\1_n)_{\ov{T}_{\pi}}$. 
 Here, Figure~\ref{baseT} are the images of homeomorphisms from the neighborhoods of 
$T_{\pi}$ and $\ov{T}_{\pi}$ to the 3-ball. 
Although $V_{\pi}$  and $V_{\pi}^{-1}$ are not unique up to ambient isotopy, 
by Lemmas~\ref{cc} and \ref{twist}, it is unique up to $C_{k}$-equivalence. 
So, for any $\pi\in\mathcal{F}_k$,  we may choose $V_{\pi}$ and $V_{\pi}^{-1}$ 
uniquely up to  $C_{k}$-equivalence. 
In particular, we may choose $V_{\pi}$ so that $\mathrm{cl}(V_{\pi})$ is the Milnor link 
$M_{\pi}=K_{\pi(1)}\cup\cdots\cup K_{\pi(k)}$ as illustrated in Figure~\ref{milnorlink2} 
(cf. Figure~\ref{milnor-tangle}). 

For $\pi\in\mathcal{F}_k$, set 
\[\mu_{\pi}(L)=\mu_L(\pi(1)\pi(2)...\pi(k))\]
By Lemma~\ref{Milnorlink}, we have the following lemma.  

\begin{lem}\label{Milnorbase1} 
For any $\pi,\pi'\in\mathcal{F}_k$, 
\[\mu_{\pi}(V_{\pi'})=
\left\{\begin{array}{ll}
1&\text{ if $\pi=\pi'$,}\\
0&\text{ if $\pi\neq\pi'$,}
\end{array}
\right.\]
and the Milnor invariants of $V_{\pi'}$ of length $\leq k-1$ vanish. 
\end{lem} 

\begin{lem}\label{base1} 
Let $T$ be a simple $C_{k-1}^d$-tree (resp. $C_{n-1}^a$-tree) for an $n$-component string link $L$. 
Then $L_T$ is $C_{k}$-equivalent (resp. $C_{n}^a$-equivalent) to $L* L'$, where 
\[L'=\prod_{\pi\in\mathcal{F}_k} V_{\pi}^{\mu_{\pi}(L_T)-
\mu_{\pi}(L)}.\]
\end{lem} 

\begin{proof}
Suppose that  $T$ is a simple $C_{k-1}^d$-tree. 
By Lemma~\ref{cc}, $L_T$ is $C_{k}$-equivalent to $L*(\1_n)_{T'}$, 
where $T'$ is a simple $C_{k-1}^d$-tree. Set $\mathrm{index}(T')=
\{i_1,...,i_k\}~(i_j<i_{j+1})$. 
Consider induction on the length of the path connecting 
the two leaves grasping $i_{k-1}$th and $i_{k}$th components of $\1_n$, and 
apply Lemma~\ref{ihx}, we have that  
$(\1_n)_{T'}$ is $C_{k}$-equivalent to 
a string link which is obtained from $\1_n$ by surgery along 
simple $C_{k-1}^d$-trees whose ends grasp $i_{k-1}$th and $i_k$th components of $\1_n$.  
By Lemmas~\ref{cc}, \ref{slide} and \ref{twist}, we have that 
\[(\1_n)_{T'}\stackrel{C_{k}}{\sim}L''=\prod_{\pi\in\mathcal{F}_k} V_{\pi}^{x_{\pi}}.\]
By Lemmas~\ref{milnCk}, \ref{add} and \ref{Milnorbase1}, 
\[\begin{array}{rcl}
\mu_{\pi'}(L_{T})=\mu_{\pi'}(L* L'')&=&\mu_{\pi'}(L)+\mu_{\pi'}(L'')\\
&=&\displaystyle  \mu_{\pi'}(L)+\sum_{\pi\in\mathcal{F}_k}x_{\pi}\mu_{\pi'}(V_{\pi})=
\mu_{\pi'}(L)+x_{\pi'}
\end{array}.\]

If $T$ is a simple $C_{n-1}^a$-tree, the arguments similar to that in the above 
can be applied. And we have the conclusion.  
\end{proof}

The following theorem gives representatives, which depend on only 
Milnor invariants, for the link-homotpy classes. 

\begin{thm}\label{SLC} 
Let $L$ be an $n$-component string link. 
Then $L$ is link-homotopic to $L_1* L_2*\cdots * L_{n-1}$, where 
\[L_i=\prod_{\pi\in\mathcal{F}_{i+1}}V_{\pi}^{x_{\pi}},~
x_{\pi}=\left\{\begin{array}{ll}
\mu_L(\pi(1)\pi(2))&\text{if $i=1$},\\ 
& \\
\mu_L(\pi(1)...\pi(i+1))-\mu_{L_1\cdots L_{i-1}}(\pi(1)...\pi(i+1))& 
\text{if $i\geq 2$}.\\
(=\mu_{L_i}(\pi(1)...\pi(i+1)))
&
\end{array}\right.
\]
\end{thm}

\begin{rem} \label{presentation} 
The presentation $L_1*L_2*\cdots*L_{n-1}$ of $L$ depends on the choice of 
order on the elements in $\mathcal{F}_i~(i=2,...,n)$. 
If we put $\mathcal{F}_2\cup\cdots\cup\mathcal{F}_{n}=\{\pi_1,...,\pi_q\}$ 
so that for $i<j$, any element in $\mathcal{F}_i$ appears before 
the elements in $\mathcal{F}_j$, then by Theorem~\ref{SLC} and Lemmas~\ref{add} and 
\ref{Milnorbase1}, 
$L$ is link-homotopic to 
$V_{\pi_1}^{x_{\pi_1}}* \cdots* V_{\pi_q}^{x_{\pi_q}}~
(x_{\pi_k}=\mu_{\pi_k}(L)-\mu_{\pi_k}(\prod_{i=1}^{k-1} V_{\pi_i}^{x_{\pi_i}}))$.   
Note that the representation is unique up to link-homotopy. 
\end{rem}

\begin{proof}
Since $C_1$-move is the crossing change, 
$L$ is $C_1$-equivalent to the trivial string link $\1_n$. 
So $L$ is obtained from $\1_n$ by surgery along 
simple $C_1$-trees.

Note that a simple $C_1$-tree is either a simple $C_1^s$-tree 
or a simple $C_1^d$-tree, and that 
$C_1^s$-equivalence preserves the value of $\mu(I)$ for any $I$ with $r(I)=1$. 
Since $L$ is $C_1^s$-equivalent to a link which is obtained from $\1_n$ by 
surgery along $C_1^d$-trees, by Lemmas~\ref{base1}, \ref{cc}, \ref{slide} and \ref{twist},   
\[L\stackrel{C_1^s+C_{2}}{\sim}\prod_{\pi\in \mathcal{F}_2}V_{\pi}^{\mu_{\pi}(L)}(=L_1).\]
A $C_k$-tree $(k\geq 2)$ is either a $C_k^{(2)}$-tree or a $C_k^d$-tree, 
and a $C_k^{(2)}$-equivalence implies $C_1^s$-equivalence (Proposition~\ref{sck}), 
and hence $(C_1^s+C_k)$-equivalence implies $(C_1^s+C_k^d)$-equivalence. 
 So $L$ is obtained form $L_1$ by surgery along 
simple $C_1^s$- and $C_2^d$-trees. 

By Lemmas~\ref{base1}, \ref{cc}, \ref{slide} and \ref{twist}, 
\[L\stackrel{C_1^s+C_{3}}{\sim}L_1*\prod_{\pi\in \mathcal{F}_3}V_{\pi}^{\mu_{\pi}(L)-
\mu_{\pi}(L_1)}(=L_1* L_2).\] 
Therefore $L$ and $L_1*  L_2$ are $(C_1^s+C_3^d)$-equivalent. 

Repeating these processes, we have that 
\[L\stackrel{C_1^s+C_{n}}{\sim}L_1*  L_2*\cdots * L_{n-1}.\] 
Since any simple $C_n$-tree for an $n$-component string link 
is a $C_n^{(2)}$-tree,  $(C_1^s+C_{n})$-equivalence
implies $C_1^s$-equivalence, i.e., link-homotopy. 
\end{proof}

By Theorem~\ref{SLC}, we have the following corollary. 

\begin{cor}\label{SLC2} 
For a natural number $k(\leq n)$, 
$n$-component string links $L$ and $L'$ are $(C_1^s+C_k)$-equivalent 
if and only if $\mu_L(I)=\mu_{L'}(I)$ for any $I$ with $r(I)=1$ and $|I|\leq k$. 
\end{cor}

\begin{proof}
The \lq only if' part follows from Lemma~\ref{milnCk}. 
Now we will prove \lq if' part. 

By Theorem~\ref{SLC}, 
$L$ and $L'$ are link-homotopic to 
$L_1*  L_2*\cdots * L_{n-1}$ and $L'_1* L'_2*\cdots * L'_{n-1}$ 
respectively. Note that both $L_i$ and $L'_i$ are $C_i$-equivalent to 
$\1_n$. So $L$ and $L'$ are $(C_1^s+C_k)$-equivalent to 
$L_1*  L_2*\cdots * L_{k-1}$ and $L'_1* L'_2*\cdots * L'_{k-1}$ 
respectively. 
Since $\mu_L(I)=\mu_{L'}(I)$ for any $I$ with $r(I)=1$ and $|I|\leq k$, 
$L_i=L_i'~(i=1,...,k-1)$.  This completes the proof. 
\end{proof} 

Theorem~\ref{LHC} follows directly from Corollary~\ref{SLC2}. 

\begin{proof}[Proof of Theorem~\ref{LHC}] 
It is enough to show \lq if' part. 
Since a $C_n$-move for an $n$-component string link is a $C_n^{(2)}$-move, 
by Proposition~\ref{sck}, $C_n$-equivalence implies $C_1^s$-equivalence. 
Hence $(C_1^s+C_n)$-equivalence implies link-homotopy. 
By Corollary~\ref{SLC2}, $L$ and $L'$ are link-homotopic. 
\end{proof} 

\begin{rem}
Let $L$ be an $n$-component link in $S^3$. 
Denote by $\mathcal{L}(L)$ the set of all $n$-component string links 
$l$ such that $\mathrm{cl}(l)=L$.  
Put $\mathcal{F}_2\cup\cdots\cup\mathcal{F}_{n}=\{\pi_1,...,\pi_q\}$ 
so that any element in $\mathcal{F}_i$ appears before 
the elements in $\mathcal{F}_j$ $(2\leq i<j\leq n)$ and fix it.  
Then, by Remark~\ref{presentation},  
each $l$ in $\mathcal{L}(L)$ is link-homotopic to 
$V_{\pi_1}^{x_{\pi_1}}* \cdots* V_{\pi_q}^{x_{\pi_q}}~
(x_{\pi_k}=\mu_{\pi_k}(l)-\mu_{\pi_k}(\prod_{i=1}^{k-1} V_{\pi_i}^{x_{\pi_i}}))$, which 
is the unique representaion up to link-homotopy. 
We define a vector $v_l$ as 
$v_l=(x_{\pi_{1}},... ,x_{\pi_{q}})$, and 
set $\mathcal{V}_L=\{ v_l ~|~ l\in \mathcal{L}(L) \}$.  
By the uniqueness of the presentation for $l$, we have the following: 
{\em Two $n$-component links $L$ and $L'$ in $S^3$ are link-homotopic  
if and only if 
$\mathcal{V}_L\cap\mathcal{V}_{L'}\neq \emptyset$.}
\end{rem}

\bigskip
{\section{Self $\Delta$-equivalence of Brunnian links } }

\begin{figure}[!h]
\includegraphics[trim=0mm 0mm 0mm 0mm, width=.8\linewidth]
{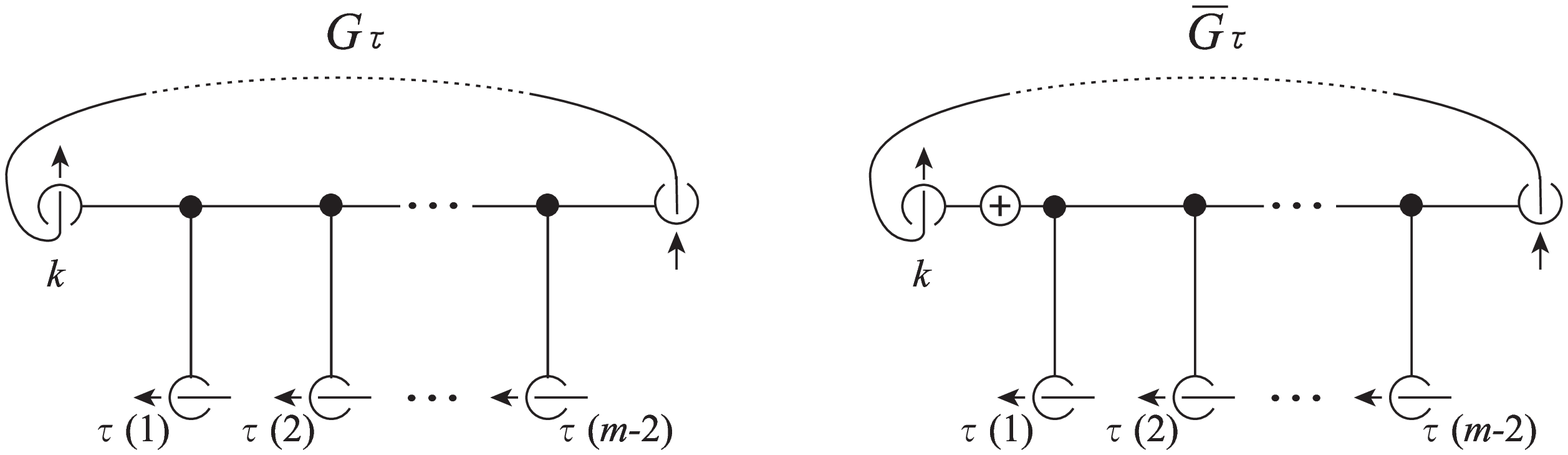}
\caption{${\tau(i)}$ and $k$ mean $\tau(i)$th and $k$th components of $\1_n$ 
respectively} \label{baseG}
\end{figure}

Let $n$ and $m$ be integers $(2\leq n< m\leq 2n)$. 
Given $k\in \n$, consider a surjection $\tau$ from $\{ 1,...,m-2 \}$ to $\n \setminus \{k \}$.  
Let $G_{\tau}$ and $\ov{G}_{\tau}$ be simple $C^a_{m-1}$-trees illustrated in Figure~\ref{baseG}, 
and  set $V_{\tau}=(\1_n)_{G_{\tau}}$ and $V^{-1}_{\tau}=(\1_n)_{\ov{G}_{\tau}}$. 

 Here, Figure~\ref{baseG} are the images of homeomorphisms from the neighborhoods of 
$G_{\tau}$ and $\ov{G}_{\tau}$ to the 3-ball. 
Although $V_{\tau}$  and $V_{\tau}^{-1}$ are not unique up to ambient isotopy, 
by Lemmas~\ref{cc} and \ref{twist}, they are unique up to $C^a_{m}$-equivalence. 
So, we may choose $V_{\tau}$ and $V_{\tau}^{-1}$ 
uniquely up to  $C^a_{m}$-equivalence.

Set
\[\mu_{\tau}(L)= \mu_L(\tau(1)...\tau(m-2)~k~k). \]  
Let $\mathcal{B}_m(k)$ be the set of all surjections $\tau$ from $\{ 1,...,m-2 \}$ to 
$\n \setminus \{k \}$ such that $|\tau^{-1}(i)|\leq 2~(i=1,...,n)$ and 
$|\tau^{-1}(j)|=1$ (if $j>k$), and 
let $\rho_m$ be a surjection from $\{ 1,...,m-2 \}$ to itself defined by $\rho_m(i)=m-1-i$.  
Note that the definition of $\mathcal{B}_m(k)$ implies that 
$\mathcal{B}_m(k)=\emptyset$ if $k<m-n$. 
So we may assume that $k\geq m-n$.

If $m\leq 2n-2$, then 
for any $\tau\in\mathcal{B}_m(k)$, $\{i~|~\tau(i)\neq\tau(m-1-i)\}\neq\emptyset$. 
We set 
\[\mathcal{P}_m(k)=\left\{\tau\in\mathcal{B}_m(k)~\left|~
\begin{array}{ll}
\tau(p)<\tau(m-1-p)~\\
\text{for}~
p=\min\{i~|~\tau(i)\neq\tau(m-1-i)\}\end{array}\right.
\right\}.\]
If $m= 2n-1$, then $k=n,n-1$ and there exists 
$\tau\in\mathcal{B}_{2n-1}(n)\cup \mathcal{B}_{2n-1}(n-1)$ such that 
$\tau(i)=\tau(2n-2-i)~(i=1,...,n-2)$ and $|\tau^{-1}(\tau(n-1))|=1$.  
For $k=n,n-1$, set 
\[\mathcal{R}_{2n-1}(k)=
\left\{\tau\in\mathcal{B}_{2n-1}(k)~\left|~
\begin{array}{ll}
\tau(i)=\tau(2n-2-i)~(i=1,...,n-2), \\
|\tau^{-1}(\tau(n-1))|=1
\end{array}\right.\right\},\]
and set
\[\mathcal{P}_{2n-1}(k)=
\left\{\tau\in \mathcal{B}_{2n-1}(k)\setminus\mathcal{R}_{2n-1}(k)
~\left|~
\begin{array}{ll}
\tau(p)<\tau(2n-2-p)~\\
\text{for}~
p=\min\{i~|~\tau(i)\neq\tau(2n-2-i)\}
\end{array}\right.\right\}.\]
Note that if $\tau\in \mathcal{R}_{2n-1}(n-1)$, then $\tau(n-1)=n$. 

If $m= 2n$, then $k=n$ and there exists $\tau\in\mathcal{B}_{2n}(n)$ such that 
$\tau(i)=\tau(2n-1-i)~(i=1,...,n-1)$.  
Set 
\[\mathcal{R}_{2n}(n)=\{\tau\in\mathcal{B}_{2n}(n)~|~\tau(i)=\tau(2n-1-i)~(i=1,...,n-1)\},\]
and set
\[\mathcal{P}_{2n}(n)=
\left\{\tau\in \mathcal{B}_{2n}(n)\setminus\mathcal{R}_{2n}(n)
~\left|~
\begin{array}{ll}
\tau(p)<\tau(2n-1-p)~\\
\text{for}~
p=\min\{i~|~\tau(i)\neq\tau(2n-1-i)\}
\end{array}\right.\right\}.\]

We note that if $\tau\in\mathcal{R}_m(k)$ then $\tau\rho_m\in\mathcal{R}_m(k)$ (i.e., $\tau$ has \lq symmetry'), 
if $\tau\in\mathcal{P}_m(k)$ then $\tau\rho_m\in\hspace*{-2ex}{/}~\mathcal{P}_m(k)$, and 
\[\mathcal{B}_m(k)={P}_m(k)\cup
\mathcal{R}_m(k)\cup\{\tau\rho_m~|~\tau\in\mathcal{P}_m(k)\}.\]
For any $\varphi\in \mathcal{B}_m(k)$, $V_{\varphi}$ is $C_{m-1}$-equivalnt 
to $\1_n$. By Lemma~\ref{milnCk}, $\mu_{V_{\varphi}}(I)=0$ for any $I$ with 
$|I|\leq m-1$. 

By the arguments similar to that in the proof of 
\cite[Proposition 5.1]{MeiY}, we have the
following lemma.

\begin{lem}\label{round} 
(1)~If $\tau\in\mathcal{P}_m(k)$, then for an $n$-string link $L$, 
\[
\mathrm{cl}(L*  V_{\tau\rho_m}) \stackrel{C_{m}^a}{\sim}
\mathrm{cl}(L*  V_{\tau}) \text{ or }
\mathrm{cl}(L*  V_{\tau}^{-1}).
\]
Moreover the $C_{m}^a$-equivalence is realized by surgery along 
simple $C_{m}^a$-trees with  
$r_j\geq r_j(G_{\tau})~(j=1,...,n)$.    \\
(2)~If $\varphi\in\mathcal{R}_{2n-1}(k)$, then for an $n$-string link $L$,
\[\mathrm{cl}(L*  V_{\varphi}) \stackrel{C_{2n-1}^a}{\sim}
\mathrm{cl}(L*  V_{\varphi}^{-1}).\]
Moreover the $C_{2n-1}^a$-equivalence is realized by surgery along 
simple $C_{m}^a$-trees with  
$r_j\geq r_j(G_{\varphi})~(j=1,...,n)$.  
\end{lem}

\begin{proof}
(1)~For the pairs of leaves grasping same components, 
by sliding the upper leaves on $\mathrm{cl}(L)$ along the orientation, 
they come to below the others. 
Lemmas~\ref{cc}, \ref{slide} and \ref{twist} complete the proof. 
The same arguments give us a proof of (2). 
\end{proof}

Now we will calculate some Milnor invariants of string links 
$V_{\varphi}$ for $\varphi\in\mathcal{R}_m(k)\cup\mathcal{P}_m(k)$.

\begin{lem}\label{Milnorbase2} 
For $\varphi\in\mathcal{P}_{m}(k)~(n+1\leq m\leq 2n,~m-n\leq k\leq n)$ and 
$\tau\in\bigcup_{l}(\mathcal{R}_{m}(l)\cup\mathcal{P}_{m}(l))$, 
\[\mu_{\tau}(V_{\varphi})=
\left\{
\begin{array}{ll}
1& \text{ if $\varphi=\tau$,}\\
0& \text{ if $\varphi\neq\tau$.}
\end{array}
\right.
\]\\
\end{lem}

\begin{proof} We take the following 4 steps to proving this lemma. 

{\bf Step 1}: 
{\em Make a new link $W_{\varphi}$ from $V_{\varphi}$ by 
taking parallels of the components of $V_{\varphi}$ so
that Lemma~\ref{parallel} can be applied.}

Let $\tau\in\mathcal{R}_m(l)\cup\mathcal{P}_m(l)~(m-n\leq l\leq n)$. 
Let $\1_m$ be the $m$-component trivial string 
link obtained from $\1_n$ by taking parallels of the components 
of $\1_n$ such that the 
$i$th component of $\1_m$ parallels to either
\[
\left\{\begin{array}{l}
\text{the $\tau(i)$th component of $\1_n$ if $i=1,...,m-2$, or}\\
\text{the $l$th component of $\1_n$ if $i=m-1,m$},
\end{array}\right.\]
and that 
$\1_m$ is contained in the tubular neighborhood $N(\1_n)$ of $\1_n$
with $G_{\varphi}\cap \1_m \subset \mathrm{int}(G_{\varphi}\cap N(\1_n))$. 
Since a surgery along $C_{m-1}$-tree preserves framings, the above correspondance 
can be naturally extended so that the
$i$th component of $(\1_m)_{G_{\varphi}}$ parallels to either 
\[
\left\{\begin{array}{ll}
\text{the $\tau(i)$th component of $(\1_n)_{G_{\varphi}}$ if $i=1,...,m-2$, or}\\
\text{the $l$th component of $(\1_n)_{G_{\varphi}}$ if $i=m-1,m$}.
\end{array}\right.\]
Set $W_{\varphi}=(\1_m)_{G_{\varphi}}$. 
By Lemma~\ref{parallel}, we have
\[\mu_{W_{\varphi}}(12... m)=\mu_{V_{\varphi}}(\tau(1)...\tau(m-2)ll)(=\mu_{\tau}(V_{\varphi})).\]

{\bf Step 2}: 
{\em By applying Lemma~\ref{split}, deform $W_{\varphi}$ up to the $C_{m}$-equivalence 
into  
\[(\1_m)_{G_{\varphi_1}}*  (\1_m)_{G_{\varphi_2}}*\cdots *(\1_m)_{G_{\varphi_s}}\] 
so that each $G_{\varphi_j}$ is a simple $C_{m-1}$-tree.}

Set $V_{j}=(\1_m)_{G_{\varphi_j}}~(j=1,...,s)$. 

{\bf Step 3}: 
{\em By applying Lemmas~\ref{milnCk} and \ref{add}, then we have 
\[\mu_{W_{\varphi}}(12...m)=\mu_{V_{1}}(12...m)+
\cdots+\mu_{V_{s}}(12...m).\] }

{\bf Step 4}: 
{\em If $G_{\varphi_j}$ is a $C_{m-1}^{(2)}$-tree, then by Proposition~\ref{sck}, 
$V_j$ is link-homotopic to trivial, hence $\mu_{V_{j}}(12...m)=0$.
Otherwise, by using Lemma~\ref{Milnorbase1}, calculate each $\mu_{V_{j}}(12...m)$.  }

If $(|\varphi^{-1}(1)|,...,|\varphi^{-1}(n)|)\neq (|\tau^{-1}(1)|,...,|\tau^{-1}(n)|)$, 
then each $G_{\varphi_j}$ is a $C_{m-1}^{(2)}$-tree. 
This implies that 
$\mu_{W_{\varphi}}(12...m)=0$.

Suppose $(|\varphi^{-1}(1)|,...,|\varphi^{-1}(n)|)= (|\tau^{-1}(1)|,...,|\tau^{-1}(n)|)$. 
Then each $G_{\varphi_j}$ is a $C_{m-1}^d$-tree. 
Since $|\varphi^{-1}(k)|=0$ and $|\tau^{-1}(i)|\geq 1~(i\neq l)$,  
we have $k=l$, i.e., $\varphi\in\mathcal{P}_m(l)$.
 
If $\varphi\neq \tau$, then neither $(\varphi(1),...,\varphi(m-1))$ nor 
$(\varphi\rho_m(1),...,\varphi\rho_m(m-1))$ is equal to 
$(\tau(1),...,\tau(m-1))$. 
By Lemma~\ref{Milnorbase1}, 
$\mu_{V_{j}}(12...m)=0$ for any $j(=1,2,...,s)$. 

If $\varphi=\tau$, then by Lemma~\ref{Milnorbase1} and by 
the fact that $\varphi\rho_m\neq\varphi$, 
there is a unique $C_{m-1}^d$-tree $G_{\varphi_u}$ in $\{G_{\varphi_1},...,G_{\varphi_s}\}$ 
such that 
\[\mu_{V_{j}}(12...m)=
\left\{\begin{array}{ll} 
1 & \text{ if $j=u$},\\
0 & \text{ if $j\neq u$}.
\end{array}
\right.\]
This completes the proof. 
\end{proof}

\begin{rem}\label{howto}
 The calculation method used in the proof of Lemma~\ref{Milnorbase2} 
can be applied for another case. 
Let $T$ be a linear, simple $C_{m-1}$-tree for $\1_n$ 
with the ends grasping $k$th component, and let 
$I=i_1...i_{m-2}kk$ be a multi-index.  
Then, $\mu_{(\1_n)_T}(I)$ can be calculated as follows.
 
{\bf Step 1}: 
{Make a new link $W=(\1_m)_{T}$ from $(\1_n)_{T}$ by 
taking parallels of the components of $(\1_n)_{T}$ so
that Lemma~\ref{parallel} can be applied.}

{\bf Step 2}: 
By applying Lemma~\ref{split}, deform $W$ up to the $C_{m}$-equivalence 
into  $(\1_m)_{T_{1}}*  (\1_m)_{T_{2}}*\cdots *(\1_m)_{T_{s}}$ 
so that each $T_{j}$ is a simple $C_{m-1}$-tree.

{\bf Step 3}: 
By applying Lemmas~\ref{milnCk} and \ref{add}, we have 
$\mu_{W}(12...m)=\mu_{(\1_m)_{T_{1}}}(12...m)+\cdots+
\mu_{(\1_m)_{T_{s}}}(12...m)$. 

{\bf Step 4}: 
If $T_{j}$ is a $C_{m-1}^{(2)}$-tree, then by Proposition~\ref{sck}, 
$\mu_{(\1_m)_{T_{j}}}(12...m)=0$.
Otherwise, by using Lemma~\ref{Milnorbase1}, calculate each 
$\mu_{(\1_m)_{T_{j}}}(12...m)$. 
\end{rem}

\begin{lem}\label{Milnorbase3} 
(1)~For any $\varphi\in\mathcal{R}_{2n-1}(k)~(k=n,n-1)$, the Milnor invariants 
of $V_{\varphi}$ of length $\leq 2n-1$ vanish.  \\
(2)~For $\varphi\in\mathcal{R}_{2n-1}(n)$ and 
$\tau\in\mathcal{R}_{2n}(n)\cup\mathcal{P}_{2n}(n)$,  
\[|\mu_{\tau}(V_{\varphi})|=
\left\{
\begin{array}{ll}
1& \text{ if $\tau\in\mathcal{R}_{2n}(n)$ and $\varphi(i)=\tau(i)~(i=1,...,n-1)$,}\\
0& \text{ otherwise.}
\end{array}
\right.
\]\\
(3)~For $\tau,~\varphi\in\mathcal{R}_{2n}(n)\cup\mathcal{P}_{2n}(n)$, 
\[\mu_{\tau}(V_{\varphi})=
\left\{
\begin{array}{ll}
1& \text{ if $\varphi=\tau\in\mathcal{P}_{2n}(n)$,}\\
2& \text{ if $\varphi=\tau\in\mathcal{R}_{2n}(n)$,}\\
0& \text{ if $\varphi\neq\tau$.}
\end{array}
\right.
\]
\end{lem}

\begin{rem} 
(1)~Note that, for any $\varphi\in\mathcal{R}_{2n-1}(n)$, there is a unique 
element $\tau\in\mathcal{R}_{2n}(n)$ such that $\varphi(i)=\tau(i)~(i=1,...,n-1)$, 
and that the correspondence induces a bijection from $\mathcal{R}_{2n-1}(n)$ to 
$\mathcal{R}_{2n}(n)$. \\ 
(2)~For $\varphi\in\mathcal{R}_{2n-1}(n-1)$ and 
$\tau\in\mathcal{R}_{2n}(n)\cup\mathcal{P}_{2n}(n)$, 
while the Milnor invariants of $V_{\varphi}$ of length $\leq 2n-1$ 
vanish, $\mu_{\tau}(V_{\varphi})$ is not easily calculated.  
And we do not need the calculations to prove Theorem~\ref{SDC}
\end{rem}

\begin{proof}
As illustrated in Figure~\ref{whitehead},
the Whitehead link, which is a link 
$C(12,12)$ 
defined in \cite[subsection 7.11]{C}, is obtained from the 2-component 
trivial link by surgery along a simple $C_2$-tree. 
We recall that, for a sequence $i_1...i_k$, \\
(i)~$C(i_1i_2,i_1i_2)$ is a Whitehead link; and\\  
(ii)~$C(i_1...i_ki_{k+1},i_1...i_ki_{k+1})=K_{i_1}\cup\cdots\cup K_{i_k+1}$ is a link obtained from  
$C(i_1...i_k,i_1...i_k)=K_{i_1}\cup\cdots\cup K_{i_{k-1}}\cup K'_{i_k}$ by replacing  
$K'_{i_k}$ with Bing doubling  $K_{i_k}\cup K_{i_{k+1}}$ of $K'_{i_k}$. 

A 4-component link obtained from the 4-component trivial 
link by surgery along the 11 {\em basic claspers} illustrated  in 
Figure~\ref{cochranlink} is ambient isotopic to 
a link illustrated in Figure~\ref{cochranlink1}, 
and  to a link obtained from the trivial link by surgery along 
a clasper with {\em boxes} as illustrated in Figure~\ref{cochranlink2}  
(for the definitions a basic clasper and a box, see \cite{H}).
Since a link illustrated in Figure~\ref{cochranlink1} is ambient isotopic to 
a link $C(1234,1234)$, the link  
$C(1234,1234)$ is $C_7$-equivalent to $\mathrm{cl}((\1_4)_G)$, where 
$G$ is a simple $C_6$-tree as illustrated in Figure~\ref{cochrantree}. 
Moreover a band sum of $C(1234,1234)$ and $\1_4$ as illustrated 
in Figure~\ref{bandsum} is $C_7$-equivalent to $(\1_4)_G$. 

Similarly, we can see that, for $\varphi\in\mathcal{R}_{2n-1}(k)~(k=n,n-1)$,  
a link $C(\alpha,\alpha)~(\alpha=\varphi(n-1)\varphi(n-2)\cdots\varphi(1)k)$ 
is $C_{2n-1}$-equivalent to either 
$\mathrm{cl}(V_{\varphi})$ or $\mathrm{cl}(V_{\varphi}^{-1})$. 
Then by Lemma~\ref{round}~(2), $C(\alpha,\alpha)$ 
is $C_{2n-1}$-equivalent to $\mathrm{cl}(V_{\varphi})$. 

In \cite[subsection 7.11]{C}, it is shown that 
$\ov{\mu}_{C(\alpha,\alpha)}(I)=0$ for any $I$ with $|I|\leq 2n-1$. 
Hence, by Lemma~\ref{milnCk}, we have the conclusion (1).

\begin{figure}[!h]
\includegraphics[trim=0mm 0mm 0mm 0mm, width=.6\linewidth]
{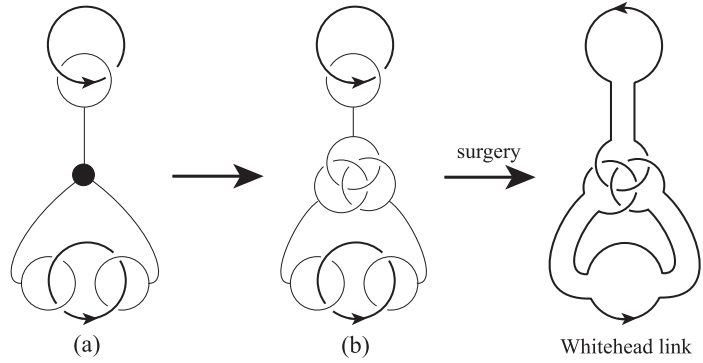}
\caption{(a): $2$-component trivial link with a simple $C_2$-clasper. 
(b): $2$-component trivial link with 3 basic claspers. } 
\label{whitehead}
\end{figure}

\begin{figure}[!h]
\includegraphics[trim=0mm 0mm 0mm 0mm, width=.5\linewidth]
{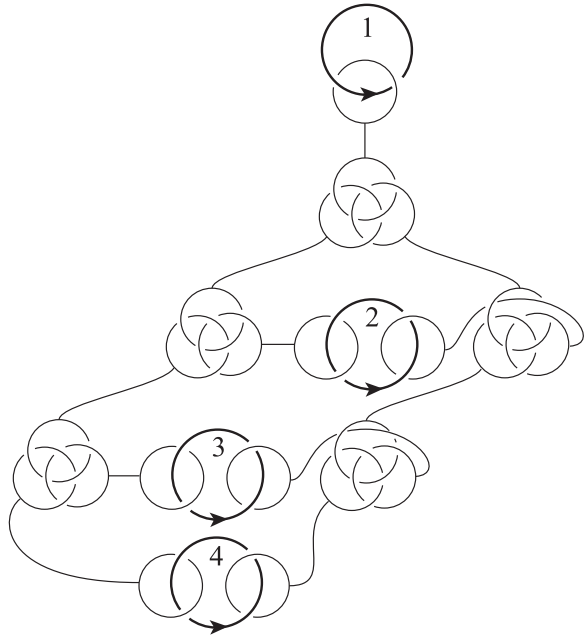}
\caption{$4$-component trivial link with 11 basic claspers. 
The numbers, 1,2,3, and 4, means the order of components. } 
\label{cochranlink}
\end{figure}

\begin{figure}[!h]
\includegraphics[trim=0mm 0mm 0mm 0mm, width=.45\linewidth]
{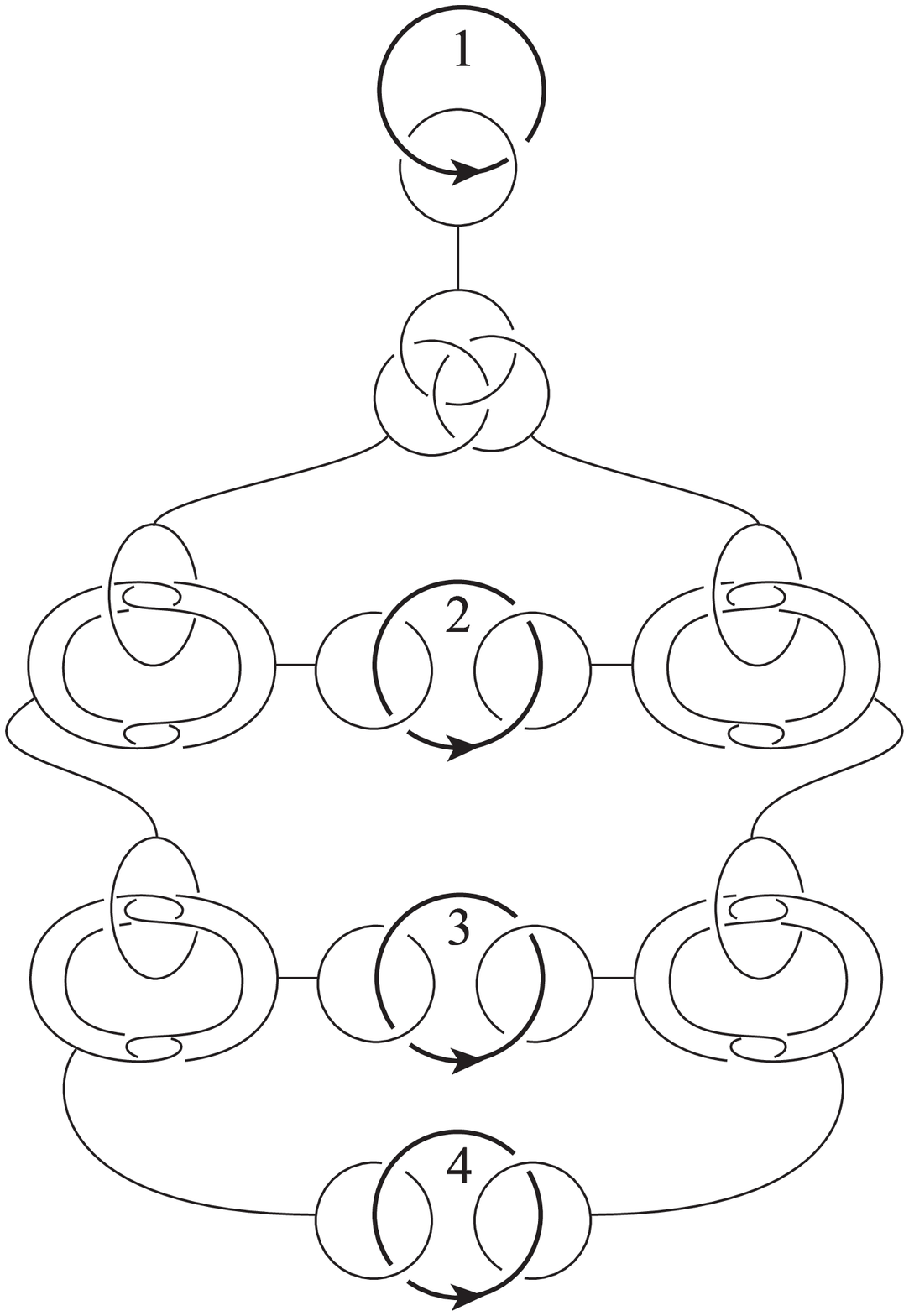}
\caption{} 
\label{cochranlink1}
\end{figure}

\begin{figure}[!h]
\includegraphics[trim=0mm 0mm 0mm 0mm, width=.4\linewidth]
{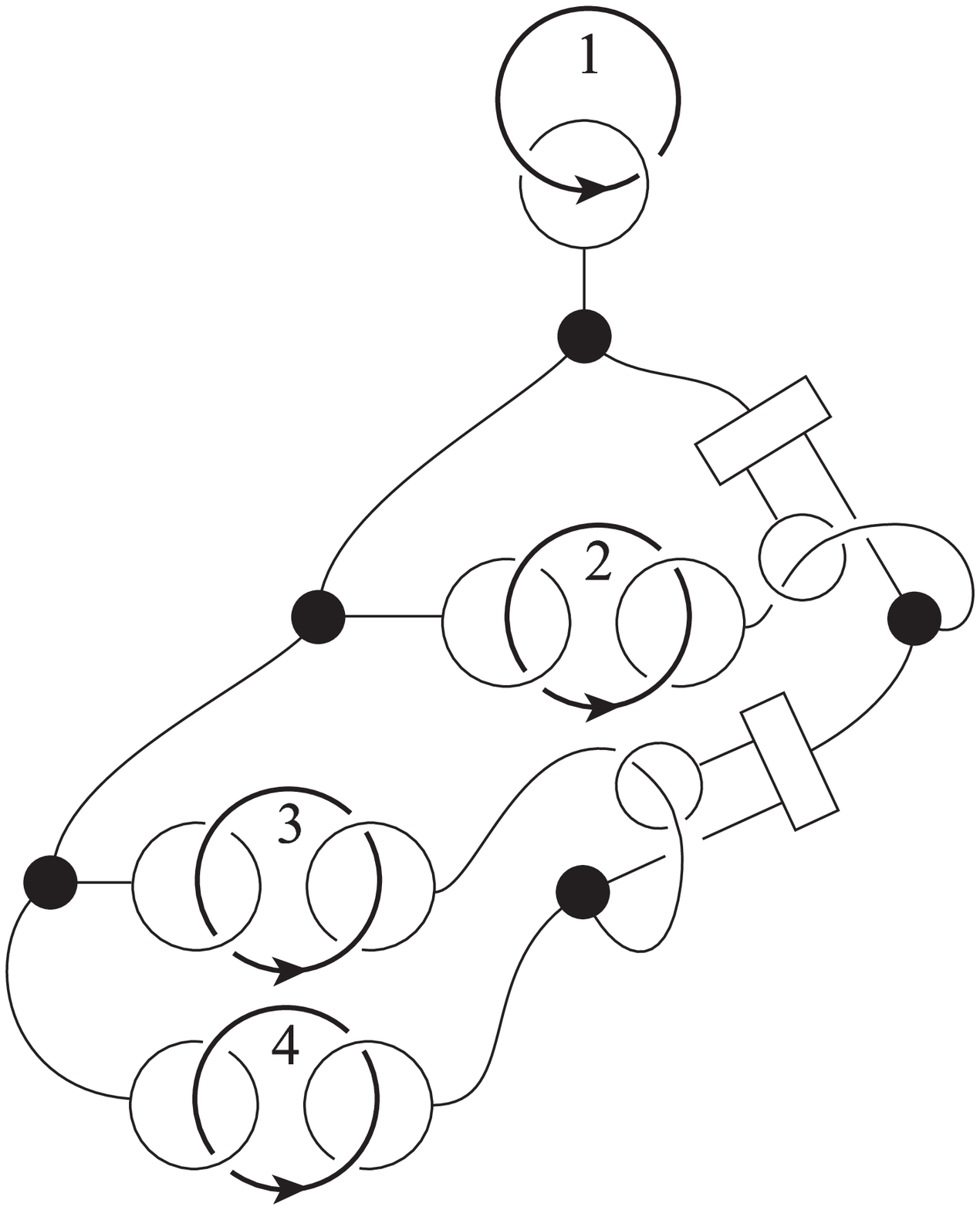}
\caption{} \label{cochranlink2}
\end{figure}

\begin{figure}[!h]
\includegraphics[trim=0mm 0mm 0mm 0mm, width=.45\linewidth]
{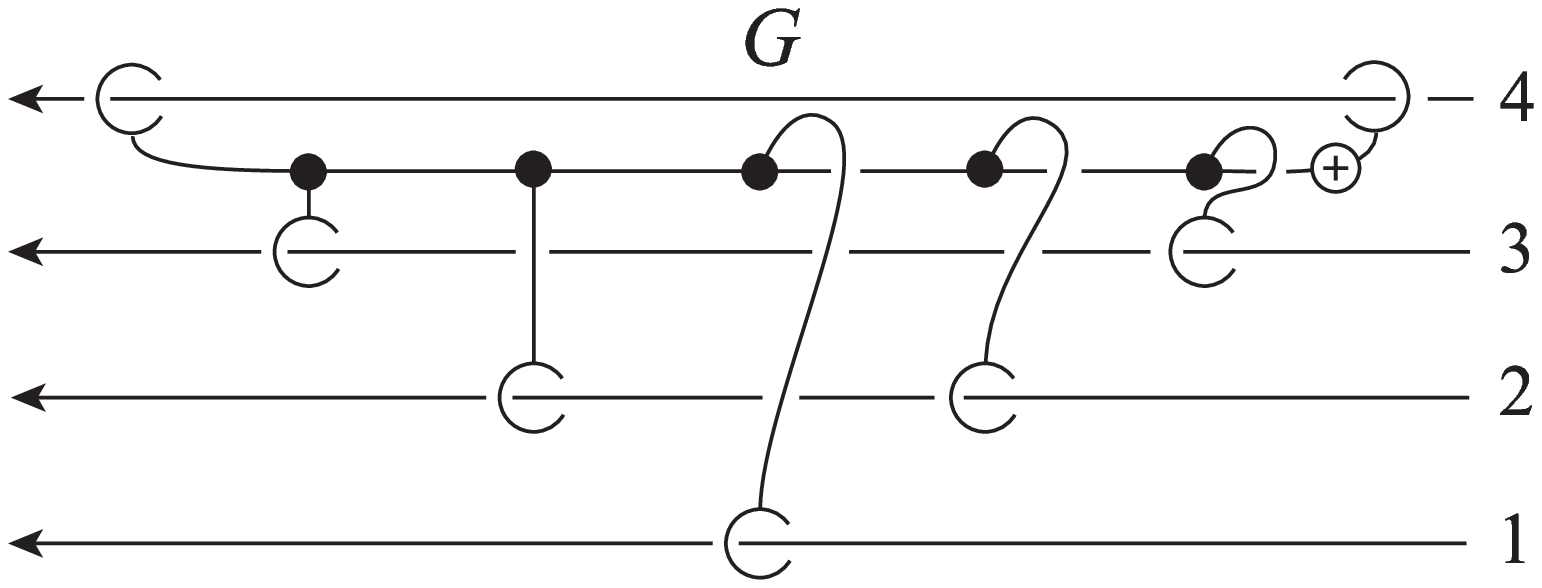}
\caption{} \label{cochrantree}
\end{figure}

\begin{figure}[!h]
\includegraphics[trim=0mm 0mm 0mm 0mm, width=.4\linewidth]
{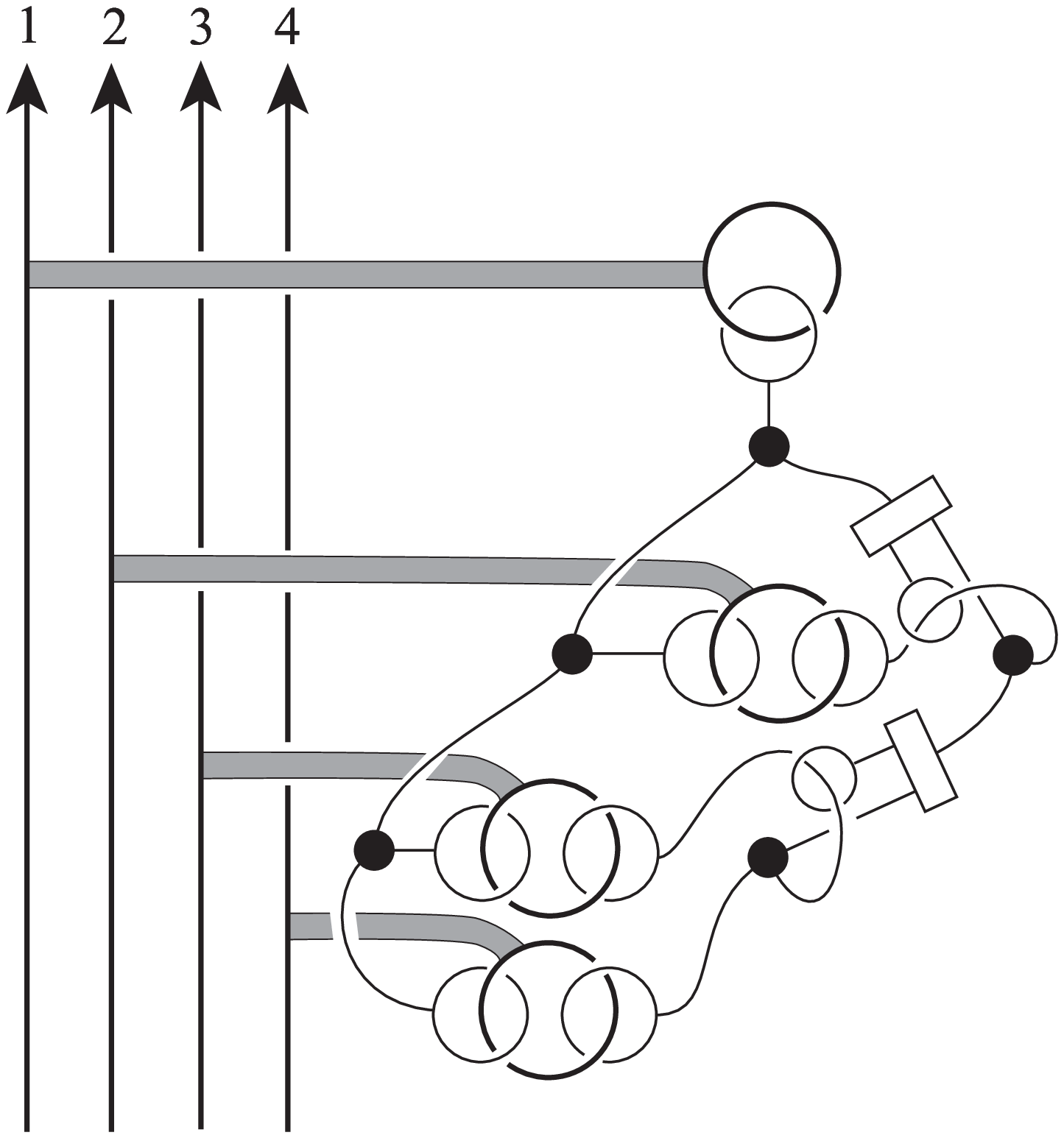}
\caption{A band sum of $\1_4$ and $C(1234,1234)$} \label{bandsum}
\end{figure}

Let $\tau\in\mathcal{R}_{2n}(n)$. 
Then it is not hard to see that 
$\mathrm{cl}(V_{\tau})$ is a link 
$L((\overline{\beta},\beta))~(~\beta=\tau(n-1)(\tau(n-2)(\cdots(\tau(1)n)\cdots)),~
\overline{\beta}=((\cdots(n\tau(1))\cdots)\tau(n-2))\tau(n-1)~)$ 
defined in \cite[subsection 7.4]{C}. 
By combining \cite[Proposition 6.5, Theorem 7.10 in p42, and Theorem 7.10 in p43]{C}, 
we have that if $\alpha=\beta$, i.e., $\varphi(i)=\tau(i)~(i=1,...,n-1)$, then 
\[|2\ov{\mu}_{C(\alpha,\alpha)}(I)|=|\ov{\mu}_{L((\overline{\beta},\beta))}(I)|\]
for any $I$ with $|I|=2n$. 
Then, (2) follows from (3). 

For $\tau,~\varphi\in\mathcal{R}_{2n}(k)\cup\mathcal{P}_{2n}(n)$, 
by following the 4 steps in Remark~\ref{howto}, we have (3).
This completes the proof. 
\end{proof}

\begin{lem}\label{base2} 
Let $L$ be an $n$-component string link and 
$m~(n+1\leq m\leq 2n)$ an integer.  
Let $T$ be a simple $C_{m-1}^a$-tree for $L$. 
Suppose that $\mu_L(I)=0$ for any $I$ with $|I|\leq m-1$ and $r(I)\leq2$, and that
$T$ is not a simple $C_{m-1}^{(3)}$-tree. Set 
$k=\max\{i~|~|T\cap(\text{$i$th component of $\1_n $})|=2\}.$ 
Then \\
(1)~$\mathrm{cl}(L_T)$ is $C_{m}^a$-equivalent to $\mathrm{cl}(L*  L')$, where 
$L'=$
\[\left\{\begin{array}{ll}
\prod_{\tau\in\mathcal{P}_m(k)} V_{\tau}^{\mu_{\tau}(L_T)-
\mu_{\tau}(L)}& \text{if $n+1\leq m\leq 2n-2$},\\
\prod_{\tau\in\mathcal{P}_{2n-1}(k)}
V_{\tau}^{\mu_{\tau}(L_T)-\mu_{\tau}(L)}
* \prod_{\varphi\in\mathcal{R}_{2n-1}(k)}
V_{\varphi}^{\varepsilon(\varphi)}& \text{if $m=2n-1$},\\
\hspace*{5cm}(\text{for some $\varepsilon(\varphi)$'s in $\{0,1\}$}) & \\
\prod_{\tau\in\mathcal{P}_{2n}(n)}
V_{\tau}^{\mu_{\tau}(L_T)-\mu_{\tau}(L)}* 
\prod_{\varphi\in\mathcal{R}_{2n}(n)}
V_{\varphi}^{(\mu_{\varphi}(L_T)-\mu_{\varphi}(L))/2}
 & \text{if $m=2n$,}
\end{array}\right.
\]
and\\
(2)~the $C_m^a$-equivalence is realized by 
surgery along simple $C_m^a$-trees with $r_k\geq 2$.  
\end{lem} 

Note that if $m\geq n+1$, then a simple $C_{m-1}^a$-tree is a simple 
$C_{m-1}^{(2)}$-tree. 
If $T$ is a simple $C_{m-1}^{(3)}$-tree, then by Proposition~\ref{sck}, 
$L_T$ and $L$ are self $\Delta$-equivalent.

\begin{proof}
By Lemma~\ref{cc}, $L_T$ is $C_m^a$-equivalent to 
$L*  (\1_n)_{T'}$, where $T'$ is a simple $C_{m-1}^a$-tree and 
not a $C_{m-1}^{(3)}$-tree with 
\[\max\{i~|~|T'\cap(\text{$i$th component of $\1_n $})|=2\}=k.\] 
By induction on the length of the path 
connecting two leaves which grasp $k$th component with 
applying Lemma~\ref{ihx}, we have that $(\1_n)_{T'}$ is $C_m^a$-equivalent 
to a string link which is obtained from $\1_n$ by surgery along 
simple linear $C_{m-1}^a$-trees whose ends grasp $k$th component of $\1_n$. 
By Lemmas~\ref{cc}, \ref{slide}, \ref{twist}, we have that 
\[L_T\stackrel{C_m^a}{\sim}L*  \prod_{\tau\in\mathcal{B}_m(k)}V_{\tau}^{x_{\tau}}.\]
Moreover, by Lemma~\ref{round}, 
$\mathrm{cl}(L*  \prod_{\tau\in\mathcal{B}_m(k)}V_{\tau}^{x_{\tau}})$
is $C_m^a$-equivalent to $\mathrm{cl}(L*  L')$, where 
\[L'=\left\{\begin{array}{ll}
\prod_{\tau\in\mathcal{P}_m(k)}V_{\tau}^{y_{\tau}} &
\text{ if $n+1\leq m\leq 2n-2$},\\
\prod_{\tau\in\mathcal{P}_m(k)}V_{\tau}^{y_{\tau}}
* \prod_{\varphi\in\mathcal{R}_m(k)}V_{\varphi}^{\varepsilon(\varphi)},
~(\varepsilon(\varphi)\in\{0,1\}) & \text{ if $m= 2n-1$},\\
\prod_{\tau\in\mathcal{P}_m(k)}V_{\tau}^{y_{\tau}}
* \prod_{\varphi\in\mathcal{R}_m(k)}V_{\varphi}^{z_{\varphi}}
 & \text{ if $m= 2n$}.
\end{array}
\right.\]
Note that the $C_m^a$-equivalences which are used in the above can be 
realized by surgery along simple $C_m^a$-trees with $r_k\geq 2$.
By Lemmas~\ref{add}, \ref{Milnorbase2}, and \ref{Milnorbase3}, 
for $\eta\in\mathcal{P}_m(k)\cup\mathcal{R}_m(k)$, we have that 
\[\mu_{\eta}(L*  L')=
\left\{
\begin{array}{ll}
\mu_{\eta}(L)+y_{\eta} & \text{ if $\eta\in \mathcal{P}_m(k)$},\\
\mu_{\eta}(L)+2z_{\eta} & \text{ if $m= 2n$ and $\eta\in \mathcal{R}_m(k)$}.
\end{array}
\right.\]
Since Milnor invariants of $L$ with length $\leq m-1$ and $r\leq2$ vanish, 
by Lemma~\ref{milnCk}, those of $L_T$ also vanish. 
By combining this, the fact that $\mathrm{cl}(L_T)$ and $\mathrm{cl}(L*L')$ are $C_m^a$-equivalence, 
and Lemma~\ref{milnCk}, we have that  
\[\mu_{\eta}(L_T)=\ov{\mu}_{\eta}(\mathrm{cl}(L_T))=
\ov{\mu}_{\eta}(\mathrm{cl}(L*L'))=\mu_{\eta}(L*L').\] 
This completes the proof.
\end{proof}

The following is the main result in this section.

\begin{thm}\label{BSD} 
Let $L$ be an $n$-component Brunnian link. 
If $\ov{\mu}_L(I)=0$ for any $I$ with $|I|\leq 2n-1$ and $r(I)\leq 2$, 
then $L$ is self $\Delta$-equivalent 
to the closure of $L'*L''$, where 
\[L'=\prod_{\varphi\in\mathcal{R}_{2n-1}(n)}V_{\varphi}^{\varepsilon(\varphi)},~~
L''=
\prod_{\tau\in\mathcal{R}_{2n}(n)}V_{\tau}^{(\ov{\mu}_{\tau}(L)-
\ov{\mu}_{\tau}(L'))/2}
*\prod_{\eta\in\mathcal{P}_{2n}(n)}V_{\eta}^{\ov{\mu}_{\eta}(L)}, 
\]
and 
\[\varepsilon(\varphi)=
\left\{
\begin{array}{ll}
1& \text{if $\ov{\mu}_{\tau}(L)$ is odd for $\tau\in\mathcal{R}_{2n}(n)$ 
with $\tau(i)=\varphi(i)~(i=1,...,n-1)$}\\
0& \text{if $\ov{\mu}_{\tau}(L)$ is even for $\tau\in\mathcal{R}_{2n}(n)$ 
with $\tau(i)=\varphi(i)~(i=1,...,n-1)$}.
\end{array}
\right.
\]
\end{thm}

Note that, 
in the theorem above, $L'*L''$ is determined by 
Milnor invariants of $L$ with length $2n$ and $r=2$. 

\begin{proof}
By Proposition~\ref{brunnian}, $L$ is obtained from the   
$n$-component trivial link $O$ by surgery along simple $C_{n-1}^a$-trees 
$T_1,...,T_l$. Hence we have 
\[L=\mathrm{cl}((\1_n)_{T_1\cup T_2\cup\cdots\cup T_l}).\] 
By Lemmas~\ref{base1}, \ref{cc}, \ref{slide} and \ref{twist}, we have that 
\[L\stackrel{C_n^a}{\sim}\mathrm{cl}(\prod_{\pi\in\mathcal{F}_n}V_{\pi}^{\ov{\mu}_{\pi}(L)}).\]

Since $\ov{\mu}_{\pi}(L)=0$ for any $\pi\in\mathcal{F}_n$, 
$L$ is $C_n^a$-equivalent to $O$, i.e.,  
$L$ is obtained from $O$ by surgery along simple $C_n^a$-trees. 
By Lemmas~\ref{base2}~(1), \ref{cc}, \ref{slide} and \ref{twist}, 
we have that 
\[L\stackrel{C_{n+1}^a}{\sim}\mathrm{cl}(\prod_{1\leq k\leq n}
(\prod_{\tau\in\mathcal{P}_{n+1}(k)}V_{\tau}^{\ov{\mu}_{\tau}(L)})).\]

Since $\ov{\mu}_{\tau}(L)=0$ for any $\tau\in\mathcal{P}_{n+1}(k)~(k=1,...,n)$, 
$L$ is $C_{n+1}^a$-equivalent to $O$. 
Note that a simple $C_m^a$-tree ($m\geq n+1$) for an $n$-component link 
is a simple $C_m^{(2)}$-tree and might be $C_m^{(3)}$-tree. 
By Lemmas~\ref{base2}~(1), \ref{cc}, \ref{slide} and \ref{twist}, 
we have that 
\[L\stackrel{C_2^s+C_{n+2}^a}{\sim}\mathrm{cl}(\prod_{2\leq k\leq n}
(\prod_{\tau\in\mathcal{P}_{n+2}(k)}V_{\tau}^{\ov{\mu}_{\tau}(L)}).\] 

Since $\ov{\mu}_{\tau}(L)=0$ for any $\tau\in\mathcal{P}_{n+2}(k)~(k=2,...,n)$, 
by repeating this step, then by Lemma~\ref{round}~(2), we have that 
\[L\stackrel{C_2^s+C_{2n-1}^a}{\sim}
\prod_{\varphi\in\mathcal{R}_{2n-1}(n)}V_{\varphi}^{\varepsilon(\varphi)}
* \prod_{\phi\in\mathcal{R}_{2n-1}(n-1)}V_{\phi}^{\varepsilon(\phi)}\]
for some $\varepsilon(\varphi)$'s and $\varepsilon(\phi)$'s  in $\{0,1\}$.

In the proof of Lemma~\ref{Milnorbase3} (1), we showed that, for 
$\phi\in\mathcal{R}_{2n-1}(n-1)$, 
\[\mathrm{cl}(V_{\phi})\stackrel{C_{2n-1}}{\sim}C(\alpha,\alpha),\]
where $\alpha=n\phi_j(n-2)\cdots\phi_j(1)(n-1)$. 
Since the Whitehead link $C(12,12)$ is deformed into a trivial link by a single self 
crossing change in the first component, 
$C(\alpha,\alpha)$ is also deformed into a trivial link by a single self 
crossing change in the $n$th component. 
So $C(\alpha,\alpha)$ is obtained from a trivial link by surgery along 
a simple $C_1^s$-tree $T$ with $r_n(T)=2$. 
By Lemma~\ref{flhb},  
$C(\alpha,\alpha)$ is obtained from a trivial link by surgery along 
simple $C_n^a$-trees with $r_n=2$. 
Since the Milnor invariants of $C(\alpha,\alpha)$ with length $\leq 2n-1$ 
vanish, by Lemmas~\ref{base2}, \ref{cc}, \ref{slide}, \ref{twist} and \ref{round}~(2), we have that 
\[C(\alpha,\alpha)\stackrel{C_2^s+C_{2n-1}^a}{\sim}
\mathrm{cl}(\prod_{\varphi\in\mathcal{R}_{2n-1}(n)}V_{\varphi}^{\varepsilon'(\varphi)})\]
for some $\varepsilon'(\varphi)$'s in $\{0,1\}$. 

Since 
\[\prod_{\varphi\in\mathcal{R}_{2n-1}(n)}V_{\varphi}^{\varepsilon(\varphi)}
* \prod_{\phi\in\mathcal{R}_{2n-1}(n-1)}V_{\phi}^{\varepsilon(\phi)}\]
is obtained from $\1_n$ by surgery along simple $C_{2n-2}^a$-trees, 
by Lemmas~\ref{cc}, \ref{slide} and \ref{round}~(2), 
$L$ is $(C_2^s+C_{2n-1})$-equivalent to the closure of 
\[L'=
\prod_{\varphi\in\mathcal{R}_{2n-1}(n)}V_{\varphi}^{\varepsilon''(\varphi)},\]
where $\varepsilon''(\varphi)$'s are integers in $\{0,1\}$. 

Note that a simple $C_{2n-1}$-tree for an $n$-component link is either $C_{2n-1}^a$-tree 
or $C_{2n-1}^{(3)}$-tree, hence by Proposition~\ref{sck}, 
$C_{2n-1}$-equivalence implies $(C_2^s+C_{2n-1}^a)$-equivalence. 
So $L$ is self $\Delta$-equivalent to a link  
obtained from $\mathrm{cl}(L')$ by surgery along simple $C_{2n-1}^a$-trees. 
By Lemmas~\ref{base2}~(1), \ref{cc}, \ref{slide} and \ref{twist},  
$L$ is $(C_2^s+C_{2n}^a)$-equivalent to $\mathrm{cl}(L'*L'')$, where 
\[L''=
\prod_{\tau\in\mathcal{R}_{2n}(n)}V_{\tau}^{(\mu_{\tau}(L)-\mu_{\tau}(L'))/2}
* \prod_{\eta\in\mathcal{P}_{2n}(n)}V_{\eta}^{\mu_{\eta}(L)}.
\]

Since a $C_{2n}$-tree is a $C_{2n}^{(3)}$-tree, by Proposition~\ref{sck}, 
$L$ is self $\Delta$-equivalent to the closure of $L'*L''$. 
By Lemmas~\ref{add}, \ref{milnCk} and \ref{Milnorbase3}, we have that 
for any $\tau\in\mathcal{R}_{2n-1}(n)$, 
\[\ov{\mu}_{\tau}(L)\equiv 
\varepsilon''(\varphi) ~(\text{mod}~ 2),\]
where 
$\varphi(i)=\tau(i)~(i=1,...,n-1)$.
This completes the proof. 
\end{proof}


By Theorem~\ref{BSD}, we have the following two 
corollaries. Corollaries \ref{BSDC} and \ref{BSDT} are 
special cases of Theorem \ref{SDC} and Corollary\ref{COR1}, respectively.
Since these corollaries are needed to show Theorem \ref{SDC}, 
we give the statements. 

\begin{cor}\label{BSDC} 
Let $L$ and $L'$ be $n$-component Brunnian links. 
Suppose that $\ov{\mu}_L(I)=\ov{\mu}_{L'}(I)=0$ for any $I$ with $|I|\leq 2n-1$ and 
$r(I)\leq 2$. 
Then $L$ and $L'$ are self $\Delta$-equivalent if and only if 
$\ov{\mu}_L(J)=\ov{\mu}_{L'}(J)$ for any $J$ with $|J|=2n$ and 
$r(J)= 2$.
\end{cor}

\begin{cor} \label{BSDT}
A Brunnian link $L$ is self $\Delta$-equivalent to a trivial link 
if and only if $\ov{\mu}_L(I)=0$ for any $I$ with $r(I)\leq 2$. 
\end{cor}
 
\bigskip
{\section{Links with Milnor invariants vanish}
}

Before proving Theorem~\ref{SDC}, we need some preparations. 

Let $L=K_1\cup\cdots\cup K_n$ be an $n$-component link and 
$b$ a band attaching a single component $K_i$ with orientation coherent, i.e., 
$b\cap L=K_i\cap b\subset \partial b$ consists of two arcs whose 
orientations from $K_i$ are opposite to those from $\partial b$. 
Then the $(n+1)$-component link $L'=(L\cup\partial b)\setminus
\mathrm{int}(b\cap K_i)$ is called a link obtained from $L$ by 
{\em fission}  (along a band $b$),  
and conversely $L$ is called a link obtained from $L'$ by {\em fusion} \cite{KSS}.
Let $L'=K_{11}\cup\cdots\cup K_{1l_1}\cup\cdots\cup K_{n1}\cup\cdots\cup K_{nl_n}$ be 
a link obtained from an $n$-component link $K_1\cup\cdots\cup K_n$ by a 
finite sequence of fission, where  
$K_{i1}\cup\cdots\cup K_{il_i}$ is obtained from $K_i~(i=1,...,n)$. 
We asign a color $c(K_{ij})$ to $K_{ij}$ as $c(K_{ij})=i$.  
In this section, for a $C_k$-tree $T$, 
we call $T$ a {\em $C_k^s$-tree} (resp. {\em$C_k^d$-tree}) if 
$|\{c(K_{ij})~|~T\cap K_{ij}\neq \emptyset\}|=1$ (resp. $=k+1$). 
A {\em $C_k^s$-move} (resp. {\em$C_k^d$-move}) is a local move defined by 
surery along simple $C_k^s$-tree (resp. $C_k^d$-tree). 

\medskip
\begin{lem}\label{fission}
If an $n$-component link is deformed into a trivial link 
by a finite sequence of fission, $C_2^s$-moves and  $C_{n-1}^d$-moves, 
then $L$ is self $\Delta$-equivalent to a Brunnian link. 
\end{lem}

\begin{proof}
Note that $L$ is obtained from a trivial link by a finite sequence of 
fusion, $C_2^s$-moves and $C_{n-1}^d$-moves. 
By the arguments similar to that in the proof of \cite[Proposition 3.22]{H}, 
we may assume that the bands of fusion, $C_2^s$-trees and $C_{n-1}^d$-trees 
are mutually disjoint. So there exist an $n$-component ribbon link $L_0$ and a  
disjoint union $F\cup F'$ of simple $C_2^s$-trees and $C_{n-1}^d$-trees 
such that $L={L_0}_{F\cup F'}$, where $F$ (resp. $F'$) is a disjoint union 
of $C_2^s$-trees (resp. $C_{n-1}^d$-trees). 
Since ribbon links are self $\Delta$-equivalent to a trivial link 
\cite{Shi}, $L_0$ is self $\Delta$-equivalent to the $n$-component 
trivial link $O$. Hence 
\[L\stackrel{C_2^s+C_{n-1}^d}{\sim}O.\]
This implies that $L$ is self $\Delta$-equivalent to a 
link obtained from $O$ by surgery along simple $C_{n-1}^d$-trees. 
Since a $C_{n-1}^d$-tree is a $C_{n-1}^a$-tree, 
by Proposition~\ref{brunnian}, we have the conclusion. 
\end{proof}

\begin{thm}\label{selfDB}
Let $L$ be an $n$-component link such that 
$\ov{\mu}_L(I)=0$ for any $I$ with $|I|\leq 2n-2$ and $r(I)\leq 2$. 
Then $L$ is self $\Delta$-equivalent to a Brunnian link.
\end{thm}

\begin{proof}
By Lemma~\ref{fission}, it is enough to show that $L$ is 
deformed into a trivial link by a finite sequence of fission, 
$C_2^s$-moves and  $C_{n-1}^d$-moves. 

Since any knot is $\Delta$-equivalent to be trivial \cite{MN}, 
we may assume that every component of $L$ is trivial. 

Suppose that any $k$-component sublink of $L$ is Brunnian 
$(2\leq k\leq n-1)$. Since surgery along a simple $C_{k-1}^d$-tree with 
index $\{i_1,...,i_k\}$ does not change the link type of 
$K_{j_1}\cup\cdots\cup K_{j_k}$ for 
$\{j_1,...,j_k\}\neq \{i_1,...,i_k\}$, by Proposition~\ref{brunnian}, 
$L$ is $C_{k-1}^d$-equivalence to a link $L'$ whose the $k$-component 
subinks are trivial. Let $L_0$ be an $n$-component string link with 
$\mathrm{cl}(L_0)=L'$, and set $\1_n=\gamma_1\cup\cdots\cup\gamma_n$.  
Let $\{S_1,...,S_m\}~(m=
\footnotesize{(\!\!\!\begin{array}{c}
n\\[-5pt]
k
\end{array}\!\!\!)})$ 
be the set of subsets of $\{1,...,n\}$ with 
$|S_i|=k~(i=1,...,m)$. 
By Lemmas~\ref{cc} and \ref{slide}, 
\[L\stackrel{C_k}{\sim}\mathrm{cl}(L_0*L_1*\cdots*L_m),\]
where $L_i$ is a string link which is 
splitable into the unions $L_{i1}$ of 
the $j$th ($j\in S_i$) components  and the other components $L_{i0}$ 
such that $L_{i1}$ is obtained from 
$\bigcup_{j\in S_i}\gamma_j$ by surgery along simple $C_{k-1}^d$-trees with 
indices $S_i$ and $L_{i0}=\1_n-\bigcup_{j\in S_i}\gamma_j$,
and the $C_k$-equivalence is realized by surgery 
along simple $C_k$-trees with $|\mathrm{index}|\geq k$.
This implies that $L$ is $C_k^d$-equivalent to a link 
obtained from $\mathrm{cl}(L_0*L_1*\cdots*L_m)$ by surgery along 
simple $C_k$-trees with $|\mathrm{index}|=k$.

By Lemmas~\ref{cc} and \ref{slide}, 
\[L\stackrel{C_k^d+C_{k+1}}{\sim}\mathrm{cl}(L_0*L_1^1*\cdots*L_m^1),\]
where $L_i^1$ is a string link which is 
splitable into the unions $L_{i1}^1$ of 
the $j$th ($j\in S_i$) components  and $L_{i0}$ such that 
$L_{i1}^1$ is obtained from $\bigcup_{j\in S_i}\gamma_j$ by surgery along 
simple $C_{k-1}^d$, $C_k$-trees with indices $S_i~(i=1,...,m)$, 
and the $C_{k+1}$-equivalence is realized by surgery 
along simple $C_{k+1}$-trees with $|\mathrm{index}|\geq k$.
By Proposition~\ref{ckindex}, surgery along a simple $C_{k+1}$-tree with 
$|\mathrm{index}|\geq k+1$ is realized by $C_k^d$-equivalence. 
Therefore, $L$ is $C_k^d$-equivalent to a link 
obtained from $\mathrm{cl}(L_0*L_1^1*\cdots*L_m^1)$ by surgery along 
simple $C_{k+1}$-trees with $|\mathrm{index}|=k$.

By Lemmas~\ref{cc} and \ref{slide}, 
\[L\stackrel{C_k^d+C_{k+2}}{\sim}\mathrm{cl}(L_0*L_1^2*\cdots*L_m^2),\]
where $L_i^2$ is 
a string link which is splitable into the unions $L_{i1}^2$ of 
the $j$th ($j\in S_i$) components  and $L_{i0}$ such that 
$L_{i1}^2$ is obtained from 
$\bigcup_{j\in S_i}\gamma_j$ by surgery along 
simple $C_{k-1}^d$, $C_k$, $C_{k+1}$-trees with indices $S_i~(i=1,...,m)$, 
and the $C_{k+2}$-equivalence is realized by surgery 
along simple $C_{k+2}$-trees with $|\mathrm{index}|\geq k$.

By repeating this procedure, we have that 
\[L\stackrel{C_k^d+C_{2k}}{\sim}\mathrm{cl}(L_0*L_1^k*\cdots*L_m^k),\]
where $L_i^k$ is a string link which is 
splitable into the unions $L_{i1}^k$ of 
the $j$th ($j\in S_i$) components  and $L_{i0}$ such that 
$L_{i1}^k$ is obtained from 
$\bigcup_{j\in S_i}\gamma_j$ by surgery along 
simple $C_{k-1}^d$, $C_k$,.., $C_{2k-1}$-trees with indices $S_i~(i=1,...,m)$.  
Note that a simple $C_{2k}$-tree is either 
a $C_{2k}$-tree with $\mathrm{index}\geq k+1$ or 
a $C_{2k}^{(3)}$-tree. 
By Lemma~\ref{sck},  
\[L\stackrel{C_k^d+C_2^s}{\sim}\mathrm{cl}(L_0*L_1^k*\cdots*L_m^k).\]
 
 So $L$ is self $\Delta$-equivalent to a link obtained from 
 $\mathrm{cl}(L_0*L_1^k*\cdots*L_m^k)$ by surgery along simple $C_k^d$-trees. 
By Lemmas~\ref{cc} and \ref{slide}, 
\[L\stackrel{C_2^s+C_{k+1}}{\sim}\mathrm{cl}(M_1*L_0*L_1^k*\cdots*L_m^k),\]
where $M_1$ is a string link obtained from $\1_n$ by surgery along 
simple $C_{k}^d$-trees, 
and the $C_{k+1}$-equivalence is realized by surgery 
along simple $C_{k+1}$-trees with $|\mathrm{index}|\geq k+1$.
 
By repeating this step, we have that
\[L\stackrel{C_2^s+C_{2n}}{\sim}\mathrm{cl}(M_{2n-k}*L_0*L_1^k*\cdots*L_m^k),\]
where $M_{2n-k}$ is a string link obtained from $\1_n$ by surgery along 
simple $C_{k}^d$, $C_{k+1}$,...,$C_{2n-1}$-trees with $|\mathrm{index}|\geq k+1$. 
Since a simple $C_{2n}$-tree is a $C_{2n}^{(3)}$-tree, by Lemma~\ref{sck}, 
\[L\stackrel{C_2^s}{\sim}\mathrm{cl}(M_{2n-k}*L_0*L_1^k*\cdots*L_m^k).\]
 
Note that $\mathrm{cl}(M_{2n-k}*L_0*L_1^k*\cdots*L_m^k)$ is deformed into 
a split sum of 
$\mathrm{cl}(M_{2n-k}*L_0)$,$\mathrm{cl}(L_{11}^k)$,...,$\mathrm{cl}(L_{m1}^k)$ 
by a finite sequence of fission. 
Since for any index $J$ in $S_i$ with $|J|\leq 2k (\leq 2n-2)$ and $r(J)\leq 2$, 
$\mu_{L_{i1}^k}(J)=\ov{\mu}_L(J)=0$, by Corollary~\ref{BSDT}, $\mathrm{cl}(L_{i1}^k)$ 
is self $\Delta$-equivalent to a trivial link. 

Hence $L$ is deformed into a split sum of $\mathrm{cl}(M_{2n-k}*L_0)$ and 
a trivial link by a finite sequence of fission and $C_2^s$-moves. 
Note that any $(k+1)$-component sublink of $\mathrm{cl}(M_{2n-k}*L_0)$ is 
Brunnian. 

By the induction, we have that 
$L$ is deformed into a split sum of an $n$-component Brunnian link $B$ and 
a trivial link by a finite sequence of fission and self $\Delta$-moves. 
By Proposition~\ref{brunnian}, $B$ is $C_{n-1}^d$-equivalent to a trivial link. 
This completes the proof. 
\end{proof}

By combining Corollarie~\ref{BSDC} and Theorem~\ref{selfDB}, 
we can prove Theorem~\ref{SDC}.

\begin{proof}[Proof of Theorem~\ref{SDC}]
Let $L$ be an $n$-component link with $\ov{\mu}_L(I)=0$ 
for any $I$ with $|I|\leq 2n-1$ and $r(I)\leq 2$. 
By Theorem~\ref{selfDB}, 
$L$ is self $\Delta$-equivalent to a Brunnian link $B$. 
Since $\ov{\mu}_{B}(I)=\ov{\mu}_{L}(I)=0$ for any 
$I$ with $|I|\leq 2n-1$ and $r(I)\leq 2$, by Corollary~\ref{BSDC}, 
$B$ is determined by Milnor invariants with length $2n$ and 
$r=2$. This completes the proof. 
\end{proof}

The following theorem characterizes $n$-component links whose Milnor 
invariants of length $\leq 2n-1$ and $r\leq 2$ vanish. 

\begin{thm}\label{char}
For an $n$-component link $L$, 
$\ov{\mu}_L(I)=0$ for any $I$ with $|I|\leq 2n-1$ and 
$r(I)\leq 2$ if and only if, for each $i\in\{1,...,n\}$, there is  
a Brunnian link $L_i$ such that 
$L_i$ is self $\Delta$-equivalent to $L$ and  
the $i$th component $K$ of $L_i$ is null-homotopic in 
$S^3\setminus(L_i-K)$.
\end{thm}

\begin{proof}
For the \lq only if' part, it is enough to consider the case when $i=n$. 
By Theorem~\ref{selfDB}, $L$ is self 
$\Delta$-equivalent to a Brunnian link. 
By Theorem~\ref{BSD}, the Brunnian link is self $\Delta$-equivalent to
the closure $L_n$ of a product of some $V_{\varphi}$'s 
($\varphi\in{\mathcal R}_{2n-1}(n)\cup{\mathcal R}_{2n}(n)
\cup{\mathcal P}_{2n}(n)$). Note that, 
for $\varphi\in{\mathcal R}_{2n-1}(n)$ 
(resp. $\varphi\in{\mathcal R}_{2n}(n)\cup{\mathcal P}_{2n}(n)$)
 $V_{\varphi}$ is $C_{2n-1}^{(2)}$-equivalent 
(resp. $C_{2n}^{(2)}$-equivalent) to $\1_n$ and the $C_{2n-1}^{(2)}$-equivalence 
(resp. $C_{2n}^{(2)}$-equivalence) is realized by surgery along 
simple $C_{2n-1}^{(2)}$-trees (resp. $C_{2n}^{(2)}$-trees) with $r_n=2$. 
By Proposition~\ref{sck}, $L_n$ is self $C_1$-equivalent to a trivia link and 
the self $C_1$-equivalence is realized by surgery along simple $C_1^s$-trees 
with $r_n=2$. Hence the $n$th component $K$ of $L_n$ is 
null-homotopic in $S^3\setminus(L_n-K)$.  

Now we will show the \lq if' part. 
Let $I$ be an index with $|I|\leq 2n-1$ and $r(I)\leq 2$. 
Since $L$ is self $\Delta$-equivalent to a Brunnian link, 
if $I$ does not contain an integer in $\{1,...,n\}$, then $\ov{\mu}_L(I)=0$. 
So we may suppose that $I$ contains any integer in $\{1,...,n\}$. 
The condition $|I|\leq 2n-1$ implies that there is an integer $i$ such that 
$i$ appears in $I$  once.  
Let $L_i$ be a Brunnian link such that 
$L_i$ is self $\Delta$-equivalent to $L$ and  
the $i$th component $K$ of $L_i$ is null-homotopic in $S^3\setminus(L_i-K)$. 
This implies that $\ov{\mu}_{L_i}(Ji)=0$ for any index $J$ in $\{1,...,n\}\setminus\{i\}$. 
Since $\ov{\mu}$ has \lq cyclic symmetry' (\cite[Theorem 8]{Milnor2}), 
$\ov{\mu}_{L_i}(I)=0$.  This completes the proof. 
\end{proof}

\begin{example}\label{example}
Let $V$ be a string link illustrated in Figure~\ref{example-figure} and 
$L$ be the closure of $V$. 
By Proposition~\ref{sck}, for $i~(i=2,3)$, $V$ is self $C_1$-equivalent 
to $\1_3$ and the self $C_1$-equivalence is realized by surgery along 
$C_1^s$-trees with indices $\{i\}$. 
Hence the $i$th component $K_i$ of $L$ is 
null-homotopic in $S^3\setminus(L-K_i)~(i=2,3)$. 
Suppose that the 1st component $K_1$ is null-homotopic in $S^3\setminus(L-K_1)$. 
Then, by Theorem~\ref{char}, 
$\ov{\mu}_L(I)=0$ for any $I$ with $|I|\leq 5$ and $r(I)\leq 2$. 
By Lemma~\ref{Milnorbase3}~(3), $\mu_V(12233)=1$. 
Hence  $\ov{\mu}_L(12233)=1$. This is a contradiction. 

\begin{figure}[!h]
\includegraphics[trim=0mm 0mm 0mm 0mm, width=.2\linewidth]
{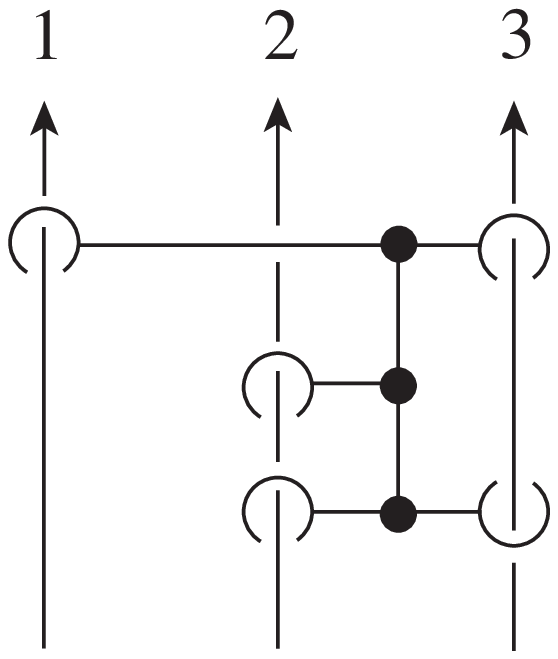}
\caption{} \label{example-figure}
\end{figure}
\end{example}

\end{document}